\documentclass[10pt,a4paper]{amsart}
\usepackage[dvips]{epsfig}
\usepackage{graphics}
\usepackage{latexsym}
\usepackage{verbatim}
\usepackage{amsmath}
\usepackage{amsthm}
\usepackage{amssymb}
\usepackage{MnSymbol}
\usepackage [table]{xcolor}
\newcommand\NN{\hbox{I\kern-.2em\hbox{N}}}
\newcommand\RR{\mathbb{R}}
\newcommand\ZZ{{{\rm Z}\kern-.28em{\rm Z}}}
\newcommand\Gradx{ \nabla_{\mathbf{x}}}
\newcommand\Gradv{ \nabla_{\mathbf{v}}}
\newcommand\Div{ \textrm{div}}
\newcommand\Rot{ \textrm{curl}}
\newcommand\jj{\mathbf{J}}
\newcommand\xx{ \mathbf{x} }
\newcommand\vv{ \mathbf{v} }
\def\signff{\bigskip\bigskip\hspace{80mm}
\vbox{{\sc Francis Filbet\par\vspace{3mm}
Universit\'e de Lyon,\par
UL1, INSAL, ECL, CNRS \par 
UMR5208, Institut Camille Jordan,\par
43 boulevard 11 novembre 1918,\par
F-69622 Villeurbanne cedex,  FRANCE
\par\vspace{3mm}e-mail:} filbet@math.univ-lyon1.fr }}

\def\signcy{\bigskip\bigskip\hspace{80mm}
\vbox{{\sc Chang Yang \par\vspace{3mm}
Universit\'e de Lyon,\par
UMR5208, Institut Camille Jordan,\par
43 boulevard 11 novembre 1918,\par
F-69622 Villeurbanne cedex,  FRANCE
\par\vspace{3mm}e-mail:} yang@math.univ-lyon1.fr }}

\title[Numerical methods based on HWENO reconstruction for Vlasov equations]{Conservative and non-conservative methods based on Hermite weighted essentially-non-oscillatory  reconstruction for Vlasov equations}\thanks{The authors are  partially supported by the European Research Council ERC Starting Grant 2009,  project 239983-\textit{NuSiKiMo}}

\author{Chang Yang and Francis Filbet}

\begin{document}

\maketitle

\begin{abstract}
We introduce a WENO reconstruction based on Hermite interpolation both for semi-Lagrangian and finite difference methods. This WENO reconstruction technique  allows to control spurious oscillations. We develop third and fifth order methods and apply them to non-conservative semi-Lagrangian  schemes and conservative finite difference methods. Our numerical results will be compared to the usual semi-Lagrangian method with cubic spline reconstruction  and the classical fifth order WENO finite difference scheme. These reconstructions  are observed  to be less dissipative than the usual  weighted essentially non-oscillatory  procedure. We apply these methods to transport equations in the context of plasma physics and the numerical simulation of turbulence phenomena. %On the one hand non-conservative semi-Lagrangian methods with high order reconstructions are particularly efficient and accurate in linear phase of simulations before the appearance of small structures.  However in the non-linear phase, the lack of conservations  may generate inaccurate numerical simulations. At contrast, conservative finite difference methods are more stable in non-linear phase and the WENO reconstruction avoid spurious oscillations.
\end{abstract}

\vspace{0.1cm}


%%%%%%%%%%%%%%%%%%%%%%%%%%%%%%%%%%%%%%%%%%%%%%%%%%%%%%%%%%%%%%%%%%%%%%%
%%                                                                   %%
%%                                                                   %%
%%%%%%%%%%%%%%%%%%%%%%%%%%%%%%%%%%%%%%%%%%%%%%%%%%%%%%%%%%%%%%%%%%%%%%%

\tableofcontents

%%%%%%%%%%%%%%%%%%%%%%%%%%%%%%%%%%%%%%%%%%%%%%%%%%%%%%%%%%%%%%%%%%%%%%%
%%                                                                   %%
%%                                                                   %%
%%%%%%%%%%%%%%%%%%%%%%%%%%%%%%%%%%%%%%%%%%%%%%%%%%%%%%%%%%%%%%%%%%%%%%%

\section{Introduction} 
\setcounter{equation}{0}
\label{sec:Intro}

Turbulent magnetized plasmas are encountered in a wide variety of astrophysical situations like the solar corona, accretion disks,  but also  in magnetic fusion devices such as tokamaks.
In practice, the study of such plasmas requires solving the Maxwell equations coupled to the computation of the plasma response.  Different ways are possible to compute this response: the fluid or the kinetic description.   Unfortunately the fluid approach seems to be insufficient when one wants to study the behavior of zonal flow, the interaction between waves and particles or the occurrence of turbulence in magnetized plasmas for example. Most of the time these plasmas are weakly collisional, and then they require a kinetic description represented by the Vlasov-Maxwell system.  The numerical simulation of the full Vlasov equation involves the discretization of the six-dimensional phase space $(\xx,\vv)\in\RR^3\times \RR^3$, which is still a challenging issue.  In the context of strongly magnetized plasmas however, the motion of the particles is particular since it is confined around the magnetic field lines;  the frequency of this cyclotron motion is faster than the frequencies of interest. Therefore, the physical system can be reduced to four or  five dimensions by averaging over the gyroradius of charged particles (See for a review \cite{ref1, ref2}).  

The development of accurate and stable numerical techniques for plasma turbulence (4D drift kinetic, 5D gyrokinetic and 6D kinetic models) is one of our long term objectives. 

Actually there are already a large variety of numerical methods based on direct numerical simulation techniques. The Vlasov equation is discretized  in phase space  using either semi-Lagrangian \cite{FSB, bibFS, bibQS11, bibSR}, finite element \cite{GQ13}, finite difference \cite{FF01, umeda,FF02} or discontinuous Galerkin \cite{bianca,heath12} schemes.  Most of these methods are based on a time splitting discretization which is particularly efficient for classical systems as Vlasov-Poisson or Vlasov-Maxwell systems.  In that case, the characteristic curves corresponding to the split operator are  straight lines and are solved exactly. Therefore, the numerical error is only due to the splitting in time and the phase space discretization of the distribution function. Furthermore for such time splitting schemes, the semi-Lagrangian methods on Cartesian grids coupled with Lagrange, Hermite or cubic spline interpolation techniques are conservative \cite{Michel,bibFS}. Hence, these methods are now currently used and have proved their efficiency for various applications. In this context semi-Lagrangian methods are often observed to be less  dissipative than classical finite volume or finite difference schemes. However, for more elaborated kinetic equations like the 4D drift kinetic \cite{bibGV} or 5D gyrokinetic \cite{ref5} equations,  or even the two dimensional guiding center  model \cite{bibSR}, time splitting techniques cannot necessarily be applied. Thus characteristic curves are more sophisticated and required a specific time discretization. For instance,  in \cite{bibGV, ref5} several numerical solvers  have been developed using an Eulerian formulation for gyro-kinetic models. However, spurious oscillations often appear  in the non-linear phase  when small structures occur and it is difficult to distinguish physical and numerical oscillations. Moreover, for these models semi-Lagrangian methods are no more conservative, hence the long time behavior of the numerical solution may become unsuitable.  

For this purpose, we want to develop a class of numerical methods based on the Hermite interpolation which is known to be less dissipative than Lagrange interpolation \cite{bibFS} together with a   weighted essentially non-oscillatory (WENO)  reconstruction applied to semi-Lagrangian and finite difference methods.
Actually, Hermite interpolation with WENO schemes were already studied in \cite{bibQS} in the context of discontinuous Galerkin methods with slope limiters. A system of equations for the unknown function and its first derivative  is evolved in time and used in the reconstruction. Moreover, a similar technique, called  CIP (Cubic Interpolation Propagation), has also been proposed for transport equations in  plasma physics applications \cite{bibCIP}, but the computational cost  is strongly increased since the unknown and all the derivatives are advected in phase space. In~\cite{bibFS}, a  semi-Lagrangian method with Hermite interpolation has been proposed and shown to be efficient and less dissipative than Lagrangian interpolation. In this latter case, the first derivatives are  approximated by a fourth order centered finite difference formula.

Here, we also apply a similar pseudo-Hermite reconstruction \cite{bibFS} and meanwhile introduce an appropriate WENO reconstruction to control spurious oscillation leading to nonlinear schemes. We develop third  and fifth order methods and apply them to semi-Lagrangian (non-conservative schemes) and conservative  finite difference methods. Our numerical results will be compared to the usual semi-Lagrangian method with cubic spline reconstruction \cite{bibSR} and the classical fifth order WENO finite difference scheme~\cite{bibJS}.

%%%%%%%%%% A REVOIR A LA FIN %%%%%%%
The paper is organized as follows : we first  present the Vlasov equation and related models which will be investigated numerically.  Then in Section \ref{sec:HSL}, the semi-Lagrangian method is proposed with high order Hermite interpolation with a WENO reconstruction to control spurious oscillations. In Section  \ref{sec:HFD},  conservative  finite difference schemes with Hermite WENO reconstructions are detailed. The Section \ref{sec:test} the one-dimensional free transport equation with oscillatory initial data is investigated to compare our schemes with classical ones (semi-Lagrangian with cubic spline interpolation and conservative finite difference schemes with WENO reconstruction). Then we perform numerical simulations on the simplified paraxial Vlasov-Poisson model and on the guiding center model for highly magnetized plasma in two dimension.

%%%%%%%%%%%%%%%%%%%%%%%%%%%%%%%%%%%%%%%%%%%%%%%%%%%%%%%%%%%%%%%%%%%%%%%%%%%%%%%%%%%%%%
%
%     SECTION 1
%
%%%%%%%%%%%%%%%%%%%%%%%%%%%%%%%%%%%%%%%%%%%%%%%%%%%%%%%%%%%%%%%%%%%%%%%%%%%%%%%%%%%%%%
          
\section{The Vlasov equation and related models}
\label{sec:vlasov}
\setcounter{equation}{0}

The evolution of the density of particles $f(t,\xx,\vv)$ in
the phase space $(\xx,\vv) \in \RR^d\times \RR^d$, $d=1,..,3,$ is given by
the Vlasov equation, 
\begin{equation}
\label{eq:vlasov} 
\frac{\partial f}{\partial t}\,+\,\mathbf{v}\cdot\Gradx f \,+\,\mathbf{F}(t,\mathbf{x},\mathbf{v})\cdot\Gradv f \,=\, 0,
\end{equation}
where the force field $F(t,\xx,\vv)$ is coupled with the distribution
function $f$ giving a nonlinear system. We mention the well known
Vlasov-Poisson (VP) and Vlasov-Maxwell (VM) models describing the evolution of particles under the effects of self-consistent electro-magnetic fields. We define the charge density $\rho(t,x)$ and current density $\jj (t,\xx)$ by
\begin{equation}
\label{densite}
\rho(t,\xx) = q\int_{\RR^d} f(t,\xx,\vv)d\vv,\quad \jj(t,\xx) = q\int_{\RR^d} \vv\,f(t,\xx,\vv)d\vv,
\end{equation}
where $q$ is the single charge. The force field is given for the Vlasov-Poisson model by
\begin{equation}
\label{poisson}
\mathbf{F}(t,\xx,\vv) = \frac{q}{m}\,\mathbf{E}(t,\xx),\quad \mathbf{E}(t,\xx) = -\nabla_{\xx} \phi(t,\xx), \quad - \Delta_\xx\phi = \frac{\rho}{\varepsilon_0},
\end{equation}
where $m$ represents the mass of one particle. For the Vlasov-Maxwell system, we have
\begin{equation}
\label{Fmaxwell}
\mathbf{F}(t,\xx,\vv) = \frac{q}{m}\,(\mathbf{E}(t,\xx) + \vv\land \mathbf{B}(t,\xx)\,),
\end{equation}
and $\mathbf{E}$, $\mathbf{B}$ are solution of the Maxwell equations
\begin{eqnarray}
\label{maxwell}
\left\{
\begin{array}{l}
\displaystyle{\frac{\partial \mathbf{E}}{\partial t} - c^2 \Rot \mathbf{B} = -\frac{\jj}{\varepsilon_0},}
\\
\,
\\
\displaystyle{\frac{\partial \mathbf{B}}{\partial t} + \Rot \mathbf{E} = 0,}
\\
\,
\\
\displaystyle{\Div \mathbf{E} = \frac{\rho}{\varepsilon_0}, \quad \Div \mathbf{B} = 0,}
\end{array}\right.
\end{eqnarray}
with the compatibility condition
\begin{equation}
\label{compatible}
\frac{\partial \rho}{\partial t} + \Div_\xx \jj = 0,
\end{equation}
which is verified by the Vlasov equation solution.

In the sequel we will also consider the so-called  guiding center model~\cite{bibCLM}, which has been derived to describe the evolution of the charge density in a highly magnetized plasma in the transverse plane of a tokamak. This model is described as follows
\begin{equation}
\left\{
\begin{array}{l}
 \displaystyle{\frac{\partial \rho}{\partial t}+\mathbf{U}\cdot\nabla\rho=0},\\[3.5mm]
 -\Delta\phi=\rho,
\end{array}
\right.
   \label{eq:gc}
  \end{equation}
where the velocity $\mathbf{U} = (-\partial_y\phi, \partial_x\phi)$ is divergence free.

Transport equations ~\eqref{eq:vlasov} or (\ref{eq:gc}) can be recast into an advective form
\begin{equation}
 \frac{\partial f}{\partial t}\,+\, \mathbf{A}\cdot\nabla f=0,
 \label{eq:transport_vlasov}
\end{equation}
where  $\mathbf{A} : \RR^{2d}\times\RR^+\to\RR^{2d}$. Hence, classical backward semi-Lagrangian method can be applied to solve~\eqref{eq:transport_vlasov}. Furthermore, under the assumption $\nabla \cdot\mathbf{A}=0$, equations \eqref{eq:vlasov} or \eqref{eq:gc} can also be rewritten in a conservative form
\begin{equation}
 \frac{\partial f}{\partial t}\,+\, \Div(\mathbf{A} f)=0.
 \label{eq:conserv_vlasov}
\end{equation}
for which a finite difference method can be used.

%%%%%%%%%%%%%%%%%%%%%%%%%%%%%%%%%%%%%%%%%%%%%%%%%%%%%%%%%%%%%%%%%%%%%%
%%                                                                   %%
%%                                                                   %%
%%%%%%%%%%%%%%%%%%%%%%%%%%%%%%%%%%%%%%%%%%%%%%%%%%%%%%%%%%%%%%%%%%%%%%%

\section{Hermite WENO reconstruction for semi-Lagrangian methods}
\label{sec:HSL}
\setcounter{equation}{0}
%%%%%%%%%%%%%%%%%%%%%%%%%%%%%%%%%%%%%%%%%%%%%%%%%%%%%%%%%%%%%%%%%%%%%

We introduce a high order Hermite interpolation coupled with a weight essentially non-oscillatory (HWENO) reconstruction for semi-Lagrangian methods. Actually, the semi-Lagrangian method  becomes a classical method for the numerical solution of the Vlasov equation because of its high accuracy and its small dissipation~\cite{bibKnorr,bibSR}. Moreover, it does not constraint any restriction on the time step size. Indeed, the key issue  of the semi-Lagrangian method compared to classical Eulerian schemes is that  it uses the characteristic curves corresponding to  the transport equation to update the unknown from one time step to the next one. Let us recall the main feature of the backward semi-Lagrangian method. For a given $s\in\RR^+$, the differential system 
\begin{equation*}
 \left\{
 \begin{array}{l}
 \displaystyle \frac{d\mathbf{X}}{dt}= \mathbf{A}(t,\mathbf{X}),\\[3mm]
 \mathbf{X}(s)=\xx,
 \end{array}
 \right.
\end{equation*}
is associated to the transport equation~\eqref{eq:transport_vlasov}. We denote its solution by $\mathbf{X}(t;s,\xx)$.
The backward semi-Lagrangian method is decomposed into two steps for computing the  function $f^{n+1}$ at time $t_{n+1}$ from the  function $f^{n}$ at time $t_{n}$ :
\begin{enumerate}
 \item For each mesh point $\xx_i$ of phase space, compute  $\mathbf{X}(t_n;t_{n+1},\xx_i)$, the value of the characteristic at time $t_n$ who is equal to $\xx_i$ at time $t_{n+1}$.
 \item As the function $f$ of transport equation verifies
 \begin{equation*}
  f^{n+1}(\xx_i)=f^n(\mathbf{X}(t_n;t_{n+1},\xx_i)),
 \end{equation*}
we obtain the value of $f^{n+1}(\xx_i)$ by computing $f^n(\mathbf{X}(t_n;t_{n+1},\xx_i))$ by interpolation, since $\mathbf{X}(t_n;t_{n+1},\xx_i)$ is not usually a mesh point.
\end{enumerate}

In practice, a cubic spline interpolation is often used~\cite{bibFG,bibGF}. It gives very good results, but it has the drawback of being non local which causes a higher communication overhead on parallel computers. Moreover  spurious oscillations may occur around  discontinuities.  On the other hand, the cubic  Hermite interpolation is local, and has been shown in~\cite{bibFS} to be less dissipative than Lagrange interpolation polynomial.  However, it has still spurious oscillations for discontinuous solution.

Here, we develop a third and fifth order Hermite interpolation coupled with a weighted essentially non-oscillatory procedure,  such that  it is accurate for smooth solutions and it removes  spurious oscillations around discontinuities or high frequencies which cannot be solved on a fixed mesh.

%%%%%%%%%%%%%%%%%%%%%%%%%%%%%%%%%%%%%%%%%%%%%%%%%%%%%%%%%%%%%%%%%%%%%%%%
\subsection{Third order Hermite WENO interpolation}

Consider a uniform mesh $(x_i)_{i}$ of the computational domain and assume that the values of the distribution function $(f_i)_i$ and its derivative $(f'_i)_i$ are known at the grid points. The standard cubic Hermite polynomial $H_3(x)$ on the interval $I_i=[x_i,x_{i+1}]$ can be expressed as follows :
\begin{equation}
\begin{array}{lll}
 H_3(x)&=& \displaystyle f_i+\frac{f_{i+1}-f_i}{\Delta x}(x-x_i)+\frac{(f_{i+1}-f_i)-\Delta x f'_i}{\Delta x^2}(x-x_i)(x-x_{i+1})\\[5mm]
 &&\displaystyle +\frac{\Delta x(f'_i+f'_{i+1})-2(f_{i+1}-f_i)}{\Delta x^3}(x-x_i)^2(x-x_{i+1}),
 \end{array}
\label{eq:poly_Hermite}
\end{equation}
 
The polynomial $H_3(x)$ verifies :
\begin{equation*}
\left\{
\begin{array}{ll}
H_3(x_i)=f_i,&H'_3(x_i)= f'_i,\\[3mm]
H_3(x_{i+1})=f_{i+1},&H'_3(x_{i+1})=f'_{i+1}.
\end{array}
\right. 
\end{equation*}

Moreover, we define two quadratic polynomials on  $I_i$ by 
$$
\left\{
\begin{array}{l}
\displaystyle h_l(x)\,=\,f_i+\frac{f_{i+1}-f_i}{\Delta x}(x-x_i)+\frac{(f_{i+1}-f_i)-\Delta x {f}'_i}{\Delta x^2}(x-x_i)(x-x_{i+1}),
\\[3mm]
\displaystyle  h_r(x)\,=\,f_i+\frac{f_{i+1}-f_i}{\Delta x}(x-x_i)+\frac{\Delta x {f}'_{i+1}-(f_{i+1}-f_i)}{\Delta x^2}(x-x_i)(x-x_{i+1}).
\end{array}\right.
$$
The polynomial $h_l$ verifies 
\begin{equation*}
 h_l(x_i)=f_i,\quad h_l(x_{i+1})=f_{i+1},\quad h'_l(x_i)= f'_i,
\end{equation*}
while  $h_r$ verifies 
\begin{equation*}
 h_r(x_i)=f_i,\quad h_r(x_{i+1})=f_{i+1},\quad h'_r(x_{i+1})= f'_{i+1}.
\end{equation*}

The idea of WENO reconstruction is now to apply  the cubic polynomial $H_3$ when the function $f$ is smooth, otherwise, we use the less oscillatory second order polynomial between $h_l$ or $h_r$. Thus, let us write  $H_3$ as follows
\begin{equation*}
 H_3(x)\,\,=\,\,w_l(x)\,p_l(x)\,\,+\,\,w_r(x)\,p_r(x),
\end{equation*}
where $w_l$ and $w_r$ are WENO weights depending on $x$.  When the function $f$ is smooth, we expect that 
 \begin{equation*}
  w_l(x)\approx c_l(x)=\frac{x_{i+1}-x}{\Delta x}\quad{\rm and}\quad w_r(x)\approx c_r(x)=1-c_l(x),
 \end{equation*}
so that we recover the cubic Hermite polynomial. Otherwise, we expect that
 \begin{equation}
  w_l(x)\approx 1, w_r(x)\approx 0\quad \text{or}\quad w_l(x)\approx 0, w_r(x)\approx 1
  \label{eq:weno_approx2}
 \end{equation}
according to the region where $f$ is less smooth. To determine these WENO weights, we follow the strategy given in~\cite{bibJS} and first define smoothness indicators by integration of the first and second derivatives of  $h_{l}$ and $h_{r}$ on the interval $I_i$ :
$$
\left\{
\begin{array}{l}
\displaystyle \beta_l\,=\,\int_{x_i}^{x_{i+1}}\Delta x(h_{l}')^2+\Delta x^3(h_{l}'')^2dx\,=\,(f_i-f_{i+1})^2+\frac{13}{3}((f_{i+1}-f_i)-\Delta x {f}'_i)^2,
\\
\,
\\
\displaystyle \beta_r\,=\,\int_{x_i}^{x_{i+1}}\Delta x(h_{r}')^2+\Delta x^3(h_{r}'')^2dx\,=\;(f_i-f_{i+1})^2+\frac{13}{3}(\Delta x {f}'_{i+1}-(f_{i+1}-f_i))^2.
\end{array}\right.
$$
Then we set $w_l$ and $w_r$ as
\begin{equation*}
 w_l(x)=\frac{\alpha_l(x)}{\alpha_l(x)+\alpha_r(x)}\quad{\rm and}\quad w_r(x)=1-w_l(x),
 %\label{eq:weno_weight1}
\end{equation*}
where
\begin{equation*}
 \alpha_l(x)=\frac{c_l(x)}{(\varepsilon+\beta_l)^2}\quad{\rm and}\quad \alpha_r(x)=\frac{c_r(x)}{(\varepsilon+\beta_r)^2}.
% \label{eq:weno_weight2}
\end{equation*}
where $\varepsilon=10^{-6}$ to avoid the denominator to be zero.

Observe that when the function  $f$ is smooth, the difference between $\beta_l$ and $\beta_r$ becomes small  and the weights  $w_l(x)\approx c_l(x)$ and $w_r(x)\approx c_r(x)$.  Otherwise, when the smoothness indicator $\beta_s$, $s=l,r$ blows-up, then the parameter $\alpha_s$ and the weight $w_s$ goes to zero, which yields~\eqref{eq:weno_approx2}.

Finally, let us mention that here the value of the first derivative  at the grid point $x_i$  is approximated by a fourth-order centered finite difference formula
 \begin{equation}
  f'_i=\frac{1}{12\Delta x}(8(f_{i+1}-f_{i-1})-(f_{i+2}-f_{i-2})).
  \label{eq:4order_dxf}
 \end{equation}
%%%%%%%%%%%%%%%%%%%%%%%%%%%%%%%%%%%%%%%%%%%%%%%%%%%%%%%%%%%%%%%%%%%%%%%%
\subsection{Fifth order Hermite WENO interpolation}
We can extend previous method to a fifth order Hermite WENO (HWENO5) interpolation.  In the same way, we first construct a fifth degree polynomial $H_5(x)$ on the interval $I_i$
$$
 H_5(x_j)=f_j, j=i-1,i,i+1,i+2,\quad H_5'(x_{i-1})=f'_{i-1}, H_5'(x_{i+2})=f'_{i+2}
$$
and  then three third degree polynomials $h_l(x)$, $h_c(x)$, $h_r(x)$ verifying
$$
\left\{
\begin{array}{lll} &h_l(x_j)=f_j, \quad j=i-1,i,i+1, & h_l'(x_{i-1})=f'_{i-1},
\\[3mm]
 &h_c(x_j)=f_j, \quad j=i-1,i,i+1,i+2,& \, 
\\[3mm]
 &h_r(x_j)=f_j, \quad j=i,i+1,i+2,& h_r'(x_{i+2})=f'_{i+2},
\end{array}\right.
$$
where the first derivative $ f'_{i}$ is given by a sixth order centered approximation
\begin{equation*}
 f'_{i}=\frac{1}{60}((f_{i+3} - f_{i-3})   -  9(f_{i+2} - f_{i-2})   +  45(f_{i+1} - f_{i-1})).
\end{equation*}

Then the polynomial  $H_5$ can be written as a convex combination
\begin{equation*}
   H_5(x) \,\,=\,\, w_l(x)h_l(x) + w_c(x)h_c(x) + w_r(x)h_r(x),
\end{equation*}
where $w_l(x)$, $w_c(x)$, $w_r(x)$ are WENO weights depending on $x$.
Similarly smoothness indicators are computed by integration of the first, second and third order derivatives of  $h_{l}(x)$, $h_c(x)$, $h_{r}(x)$ on the interval $I_i$ :
\begin{equation*}
 \beta_j = \int_{x_i}^{x_{i+1}}\Delta x(h_{j}')^2+\Delta x^3(h_{j}'')^2 + \Delta x^5(h_{j}''')^2 dx,\quad j=l,c,r.
\end{equation*}
Finally, the WENO weights are determined according to the smoothness indicators
$$
\left\{
\begin{array}{lll}
 \displaystyle w_l(x) = \frac{\alpha_l(x)}{\alpha_l(x) + \alpha_c(x) + \alpha_r(x) }, & \displaystyle\alpha_l(x) = \frac{c_l(x)}{(\varepsilon + \beta_l)^2}, & \displaystyle c_l(x)=\frac{(x-x_{i+2})^2}{9\Delta x^2},\\[5mm]
\displaystyle   w_c(x) = \frac{\alpha_c(x)}{\alpha_l(x) + \alpha_c(x) + \alpha_r(x) }, & \displaystyle\alpha_c(x) = \frac{c_c(x)}{(\varepsilon + \beta_c)^2}, & \displaystyle c_c(x)=1-c_l(x)-c_r(x), \\[5mm]
 \displaystyle  w_r(x) = \frac{\alpha_r(x)}{\alpha_l(x) + \alpha_c(x) + \alpha_r(x) }, & \displaystyle\alpha_r(x) = \frac{c_r(x)}{(\varepsilon + \beta_r)^2}, & \displaystyle c_r(x)=\frac{(x-x_{i-1})^2}{9\Delta x^2}.
\end{array}\right.
$$

This polynomial reconstruction allows to get fifth order accuracy  for smooth stencil and the various stencils are expected to damp oscillations when filamentation of the distribution function occurs. Finally, let us observe that this technique can be easily extended to high space dimension on Cartesian grids.
  
\section{Hermite WENO reconstruction for conservative finite difference methods}
\label{sec:HFD}
\setcounter{equation}{0}
When the velocity $\mathbf{A}$ is not constant (\ref{eq:transport_vlasov}), the semi-Lagrangian method is not conservative even when $\Div \mathbf{A}=0$, hence mass is no longer conserved and the long time behavior of the numerical solution can be wrong even for small time steps. Therefore, high order conservative methods may be more appropriate even if they are restricted by a CFL condition. An alternative is to use the finite difference formulation in the conservative form and to use the semi-Lagrangian method for the flux computation \cite{Russo,Qiu}.

In this section, we extend Hermite WENO reconstruction for computing numerical flux of finite difference method.
Suppose that $\{f_i\}_{1\leq i\leq N}$  is approximation of $f(x_i)$. 
We look for the flux $\{\hat f_{i+1/2}\}_{0\leq i\leq N}$ such that it approximates the derivative $f'(x)$ to $k$-th order accuracy :
\begin{equation*}
 \frac{ \hat f_{i+1/2} - \hat f_{i-1/2} }{\Delta x} = f'(x) + \mathcal{O} ( \Delta x^k ).
\end{equation*}
Let us define a function $g$ such that
\begin{equation}
 f(x) = \frac{1}{\Delta x} \int_{x-\Delta x/2}^{x + \Delta x/2} g(s) ds,
 \label{eq:intro_h}
\end{equation}
then clearly
\begin{equation*}
 f'(x) = \frac{1}{\Delta x}\left[ g(x + \Delta x/2) - g(x - \Delta x/2) \right].
\end{equation*}
Hence we only need
\begin{equation*}
 \hat f_{i+1/2} \approx p( x_i + \Delta x/2 ).
\end{equation*}
Let us denote by $G$ one primitive of $g$ 
\begin{equation*}
 G(x) = \int_{-\infty}^x g(s) ds,
\end{equation*}
then~\eqref{eq:intro_h} implies
\begin{equation*}
 G(x_{i+1/2} ) = \sum_{j=-\infty} ^i \int_{x_{j-1/2}} ^{x_{j+1/2}} g(s) ds = \Delta x \sum_{j=-\infty} ^i f_j \,=:\, G_{i+1/2}.
\end{equation*}
Thus, given the point values $\{f_i\}_i$, the primitive function $G(x)$ is exactly known at $x=x_{x+1/2}$.
We thus can approximate  $G(x)$ by an interpolation method. 
Therefore, 
\begin{equation}
 g(x_{i+1/2}) = \left.\frac{d G}{d x}\right|_{x=x_{i+1/2}}.
\end{equation}

Now let us interpolate the primitive function $G(x)$. Here we give the Hermite WENO scheme and outline the procedure of reconstruction only for the fifth order accuracy case.

The aim is to construct an approximation of the flux  $f^{-}_{i+1/2}$ by the Hermite polynomial of degree five together with a WENO reconstruction from point values $\{f_i\}$ :
\begin{enumerate}
 \item We construct the Hermite polynomial $H_5$ such that
$$
H_5(x_{j+1/2}) = G_{j+1/2}, j=-2,-1,0,1,\quad H'_5(x_{j+1/2}) =  G'_{j+1/2}, j=-1,0,
$$
\item We construct  cubic reconstruction polynomials $H_l(x)$, $H_c(x)$, $H_r(x)$ such that :
 $$
\left\{
\begin{array}{ll}
 H_l(x_{j+1/2}) = G_{j+1/2}, j=-2,-1,0, & H'_l(x_{i-1/2}) =  G'_{i-1/2},\\[4mm]
 H_c(x_{j+1/2}) = G_{j+1/2}, j=-2,-1,0,1,& \, \\[4mm]
 H_r(x_{j+1/2}) = G_{j+1/2}, j=-1,0,1, & H'_r(x_{i+1/2}) =  G'_{i+1/2},\\[4mm]
  \end{array}\right.
$$
where $G'_{i+1/2}$ is the sixth order centered approximation of first derivative. Let us denote by $h_l(x)$, $h_c(x)$, $h_r(x)$, $h_5(x)$  the first derivatives of $H_l(x)$, $H_c(x)$, $H_r(x)$, $H_5(x)$  respectively. By evaluating  $h_l(x)$, $h_c(x)$, $h_r(x)$, $h_5(x)$ at $x=x_{i+1/2}$, we obtain
$$
 h_5(x_{i+1/2}) \,=\, \frac{ -8f_{i-1} + 19f_{i} + 19f_{i+1}   + 3H'_{i-1/2} - 6H'_{i+1/2}}{27}
$$
and
\begin{eqnarray*}
 h_l(x_{i+1/2}) &=&  -2\,f_{i-1} \,+\, 2f_{i} \,+\, G'_{i-1/2} , \\[3mm]
 h_c(x_{i+1/2}) &=& \frac{-f_{i-1} \,+\, 5\,f_i + 2\,f_{i+1}}{6},\\[3mm]
 h_r(x_{i+1/2}) &=& \frac{ f_{i} \,+\, 5\,f_{i+1}  \,-\, 2\,G'_{i+1/2}}{4}.
\end{eqnarray*}

\item We evaluate the smoothness indicators $\beta_l$, $\beta_c$, $\beta_r$, which measure the smoothness of $h_l(x)$, $h_c(x)$, $h_r(x)$ on the cell $I_i$.
\begin{eqnarray*}
 \beta_l &=&\int_{x_i}^{x_{i+1}}  \Delta x (h'_l(x))^2 +  \Delta x^3 (h''_l(x))^2 dx\\[3mm]
 &=& \frac{1}{16}\left(835f_{i-1}^2 + 139f_i^2 + 300(H'_{i-1/2})^2 - 674f_{i-1}f_i - 996f_{i-1}H'_{i-1/2} + 396f_i H'_{i-1/2}\right),\\[3mm]
 \beta_c &=&\int_{x_i}^{x_{i+1}}  \Delta x (h'_c(x))^2 +  \Delta x^3 (h''_c(x))^2 dx\\[3mm]
 &=& \frac{1}{12}\left(13f_{i-1}^2 + 64f_i^2 + 25f_{i+1}^2 - 52f_{i-1}f_i + 26f_{i-1}f_{i+1} - 76f_i f_{i+1}\right),\\[3mm]
 \beta_r &=&\int_{x_i}^{x_{i+1}}  \Delta x (h'_r(x))^2 +  \Delta x^3 (h''_r(x))^2 dx\\[3mm]
 &=& \frac{1}{16}\left(55f_{i}^2 + 367f_{i+1}^2 + 156(H'_{i+1/2})^2 - 266f_{i}f_{i+1} + 156f_{i}H'_{i+1/2} - 468f_{i+1} H'_{i+1/2}\right).
\end{eqnarray*}

\item We compute the non-linear weights based on the smoothness indicators
$$
\left\{
\begin{array}{ll}
 \displaystyle w_l = \frac{\alpha_l}{\alpha_l + \alpha_c + \alpha_r }, &  \displaystyle\alpha_l = \frac{c_l}{(\varepsilon + \beta_l)^2}, \\[5mm]
  \displaystyle w_c = \frac{\alpha_c}{\alpha_l + \alpha_c + \alpha_r }, &  \displaystyle\alpha_c = \frac{c_c}{(\varepsilon + \beta_c)^2}, \\[5mm]
 \displaystyle  w_r = \frac{\alpha_r}{\alpha_l + \alpha_c + \alpha_r }, &  \displaystyle\alpha_r = \frac{c_r}{(\varepsilon + \beta_r)^2}, 
\end{array}\right.
$$
where the coefficients $c_l=1/9$, $c_c=4/9$, $c_r=4/9$ are chosen to get fifth order accuracy for smooth solutions and the parameter $\varepsilon=10^{-6}$ avoids the blow-up of $\alpha_k$, $k=\{l,c,r\}$.

\item The flux $f^-_{i+1/2} $ is then computed as
\begin{equation*}
 f^-_{i+1/2} \,\,=\,\; w_l \,h_l(x_{i+1/2})  \,\,+\,\, w_c\, h_c(x_{i+1/2})  \,\,+\,\, w_r \,h_r(x_{i+1/2}). 
\end{equation*}

\end{enumerate}
The reconstruction to $f^+_{i+1/2}$ is mirror symmetric with respect to $x_{i+1/2}$ of the above procedure.
%%--------------------------------------------------------------------
\section{Numerical simulation of Vlasov equation and related models}
\label{sec:test}
\setcounter{equation}{0}
We start with a very basic test on the one dimensional transport equation with constant velocity to check the order of accuracy and to compare the error amplitude of the various numerical schemes. Then we perform numerical simulations on the simplified paraxial Vlasov-Poisson model and on the guiding center model for  highly magnetized plasma in the transverse plane of a tokamak.

In this section we will compare our Hermite WENO reconstruction with the usual semi-Lagrangian method with cubic spline interpolation without splitting \cite{bibSR},  and with the classical fifth order finite difference technique \cite{bibJS} coupled with a fourth order Runge-Kutta scheme for the time discretization.    
\subsection{1D transport equation}
\label{sec:1DTransport}
We compare our Hermite WENO reconstruction with various classical methods for solving the free  transport equation
\begin{equation}
 \frac{\partial f}{\partial t}\,\,+\,\,\,  \frac{\partial f}{\partial x}\,\,=\,\,0,\quad x\in[-1,1],\quad t\geq0,
 \label{eq:test1D}
\end{equation}
with periodic boundary conditions.

 Let us first consider a smooth solution, where the  initial condition is chosen as
\begin{equation*}
 f(0,x)=\sin\left(\pi x\right),\quad x\in[-1,1].
\end{equation*} 
We present in Table~\ref{tab:1D1}, the numerical error for  different methods. On the one hand for semi-Lagrangian methods, the Hermite WENO interpolation is compared with the cubic spline interpolation. The semi-Lagrangian method is unconditionally stable, we thus choose a CFL number larger than one, {\it e.g.} CFL $=2.5$. We observe that the cubic spline and Hermite WENO reconstructions have both third order accuracy, and the numerical error has almost the same amplitude. The semi-Lagrangian method with a fifth order  Hermite WENO reconstruction has fifth order accuracy, thus it is much more accurate than the previous third order methods.

On the other hand we focus on the finite difference method and compare the Hermite WENO reconstruction  with the classical fifth order WENO reconstruction~\cite{bibJS}. We observe that these two methods have fifth order accuracy, but the Hermite WENO interpolation method is much more accurate than the usual WENO method. Furthermore, for the same order of accuracy the semi-Lagrangian method is much more precise than  the finite difference scheme, which is expected for linear problems since the error only comes from the polynomial interpolation.

\begin{table}[h]
\begin{center}
 \begin{tabular}{|c|c|c|c|c|c|c|}
  \hline
  $n_x$&\multicolumn{2}{c|}{200}  & \multicolumn{2}{c|}{400}   & \multicolumn{2}{c|}{800} \\
  \cline{2-7}
  & $\|\cdot\|_1$  & $r$  & $\|\cdot\|_1$  & $r$  & $\|\cdot\|_1$  & $r$ \\
  \hline  
  Semi-Lagrangian cubic spline&   1.03e-6   &  3.00   & 1.29e-7   & 3.00      &   1.61e-8  &  3.00     \\
  \hline
  Semi-Lagrangian cubic HWENO&    1.04e-6 &  3.03  &   1.29e-7     &  3.01 &   1.62e-8  &  3.00   \\
  \hline
  Semi-Lagrangian HWENO 5th &     9.28e-10 &  5.51  &   2.28e-11     &  5.35 &   8.63e-13  &  4.72   \\
  \hline
  Finite difference WENO 5th&    1.15e-7 &  4.99   & 3.60e-9   & 5.01   & 1.16e-10   & 4.95     \\
  \hline
  Finite difference  HWENO 5th&    6.06e-8 &  4.99  &   1.92e-9     &  4.98 &   6.55e-11  &  4.87    \\
  \hline
 \end{tabular}
\end{center}
\caption{\label{tab:1D1}1D transport equation : {\it Error  in $L_1$-norm and order of convergence $r$ for smooth solutions for semi-Lagrangian and finite difference methods. The  final time is $T_{\text{end}}=8$.}}
\end{table}

We next  consider a step function as follows
\begin{equation}
 f(0,x)=
 \left\{
 \begin{array}{ll}
  1,&\textrm{ for } -1\leq x\leq0,\\[3mm]
  0,&\text{otherwise}.
 \end{array}
 \right.
  \label{eq:discontinuous_initial1}
\end{equation}

Comparisons of the methods are now summarized in Table~\ref{tab:1D2}. We first notice that all the methods can  achieve  order of accuracy of $\frac{p}{p+1}$, where $p$ is degree of polynomial. On the one hand, it is clear that semi-Lagrangian methods are more precise than finite difference ones. More precisely for $n_x\leq 10^3$, the cubic spline interpolation is more accurate than Hermite of degree three and five coupled with the WENO reconstruction. It illustrates perfectly the robustness of the semi-Lagrangian method with cubic spline interpolation.  Nevertheless it also generates a lot of oscillations (see  Table~\ref{tab:1D2} (b)) which produce negative values of the distribution function. Furthermore the Hermite WENO5 reconstruction is less dissipative than usual WENO5 and it is more accurate. Both of them   control well spurious oscillations (see   Table~\ref{tab:1D2} (b)).

\begin{table}[h]
\begin{center}
 \begin{tabular}{|c|c|c|c|c|c|c|}
  \hline
  $n_x$&\multicolumn{2}{c|}{200}  & \multicolumn{2}{c|}{400}   & \multicolumn{2}{c|}{800} \\
  \cline{2-7}
  & $\|\cdot\|_1$  & $r$  & $\|\cdot\|_1$  & $r$  & $\|\cdot\|_1$  & $r$ \\
  \hline  
  Semi-Lagrangian cubic spline&   2.47e-2   &  0.80   & 1.43e-2   & 0.79      &   8.52e-2  &  0.75     \\
  \hline
  Semi-Lagrangian-HWENO 3rd&    3.27e-2 &  0.78  &   1.89e-2     &  0.79 &   1.09e-2  &  0.80   \\
  \hline
  Semi-Lagrangian-HWENO 5th&    2.94e-2 &  0.84  &   1.63e-2     &  0.85 &   8.99e-3  &  0.86   \\
  \hline
  Finite difference-WENO 5th&    4.50e-2 &  0.83   & 2.53e-2   & 0.83   & 1.43e-2   & 0.82     \\
  \hline
  Finite difference-HWENO 5th&    4.07e-2 &  0.83  &   2.29e-2     &  0.83 &   1.29e-2  &  0.83    \\
  \hline
 \end{tabular}
 \vskip 3mm
 (a) Error between exact solution and approximated solution
  \vskip 3mm
   \begin{tabular}{|c|c|c|c|}
  \hline
  $n_x$&200&400&800\\
  \hline
  Semi-Lagrangian  cubic spline&  5.75e-1      &  5.12e-1      & 5.19e-1  \\
  \hline
  Semi-Lagrangian-HWENO 3rd&  5.71e-4      &  5.18e-4      & 4.37e-4    \\
  \hline
  Semi-Lagrangian-HWENO 5th&  1.09e-3      &  1.42e-3      & 1.46e-3    \\
  \hline
  Finite difference-WENO 5th&  9.54e-5      &  8.43e-5  &  6.58e-5    \\
  \hline
  Finite difference-HWENO 5th&  2.30e-3      &  2.47e-3  &  1.88e-3  \\
  \hline
 \end{tabular}
  \vskip 3mm
 (b) Error of total variation
  \vskip 3mm
 \caption{ \label{tab:1D2}1D transport equation : {\it  Comparison of different methods for the linear equation~\eqref{eq:test1D} with initial data~\eqref{eq:discontinuous_initial1}. (a) Error in $L_1$ norm and $r$ is the order of accuracy (b) Error on the total variation. The  final time is $T_{\text{end}}=8$.}}
 \end{center}
\end{table}

We finally consider an oscillatory solution where the initial condition is given by~\cite{bibJS}, 
\begin{equation}
 f(0,x)=
 \left\{
 \begin{array}{ll}
  \frac{1}{6}[G(x,z-\delta)\,+\,G(x,z-\delta)\,+\,4\,G(x,z)],&\textrm{ for } -0.8\leq x\leq-0.6,\\[3mm]
  1,&\textrm{ for } -0.4\leq x\leq-0.2,\\[3mm]
  1-|10(x-0.1)|, &\textrm{ for } 0\leq x\leq0.2,\\[3mm]
   \frac{1}{6}[F(x,z-\delta)\,+\,F(x,z-\delta)\,+\,4\,F(x,z)],&\textrm{ for } 0.4\leq x\leq0.6,\\[3mm]
  0,&\text{otherwise}.
 \end{array}
 \right.
  \label{eq:discontinuous_initial2}
\end{equation}
where $G(x,z)=\exp(-\beta(x-z)^2)$, $F(x,a)=\{\max((1-\alpha^2(x-a)^2)^{1/2},0)\}$ with  $\alpha=0.5$, $z=-0.7$, $\delta=0.005$, $\alpha=10$ and $\beta=(\log2)/36\delta^2$.

\begin{figure}
\begin{center}
 \begin{tabular}{cc}
  \includegraphics[width=6.5cm]{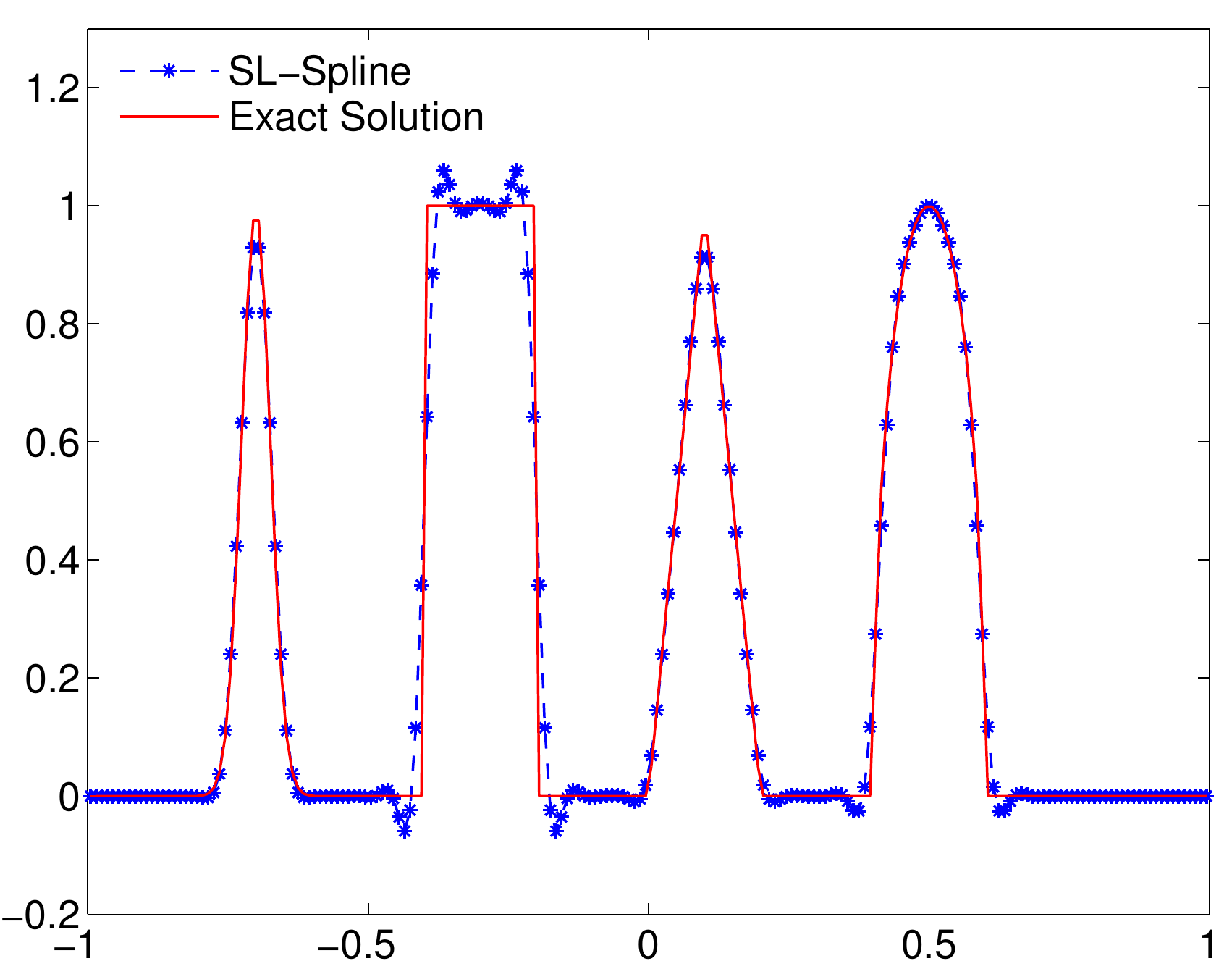} &    \includegraphics[width=6.5cm]{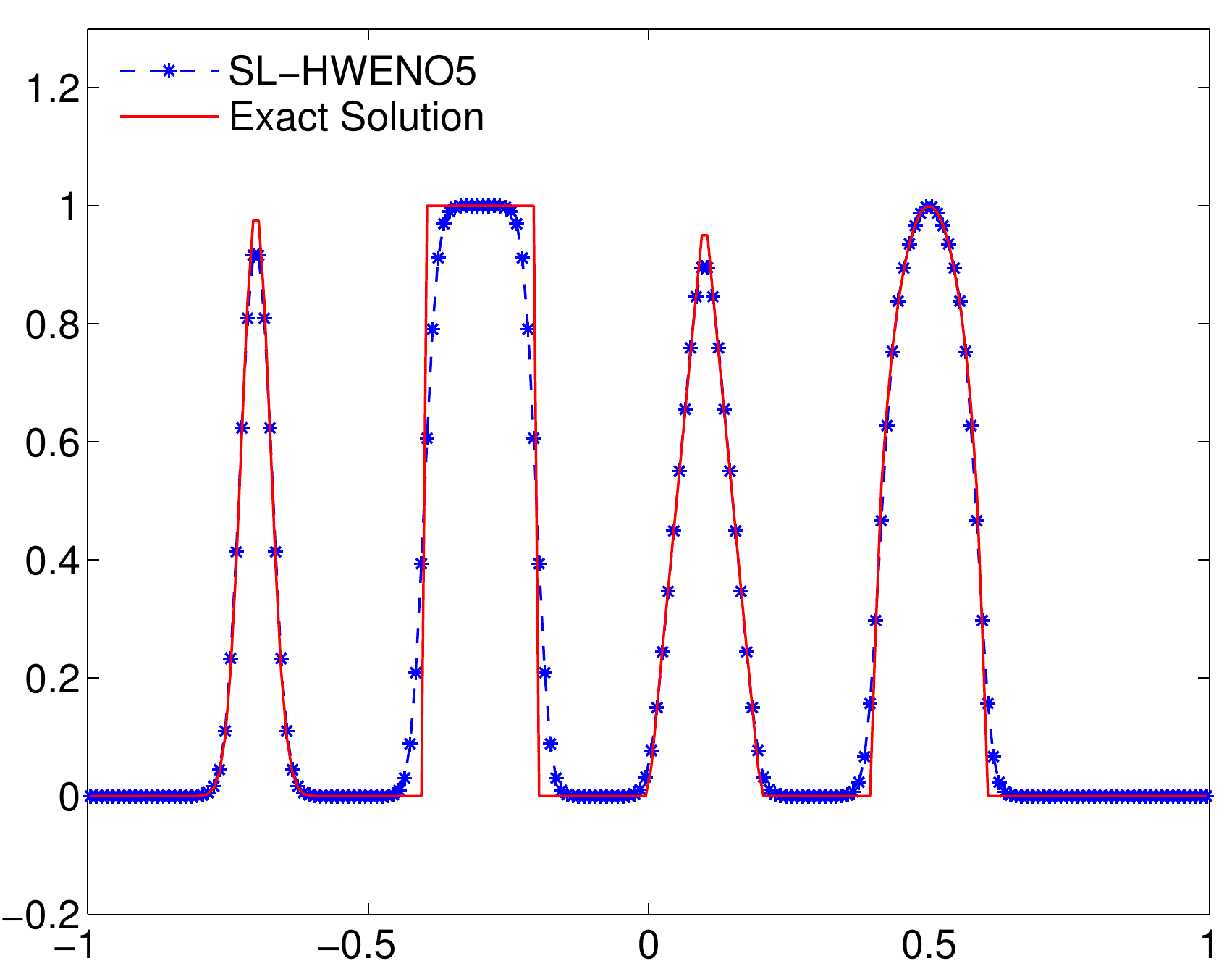} \\
  (a) Semi-Lagrangian cubic spline                 &    (b) Semi-Lagrangian-HWENO 5th \\
  \includegraphics[width=6.5cm]{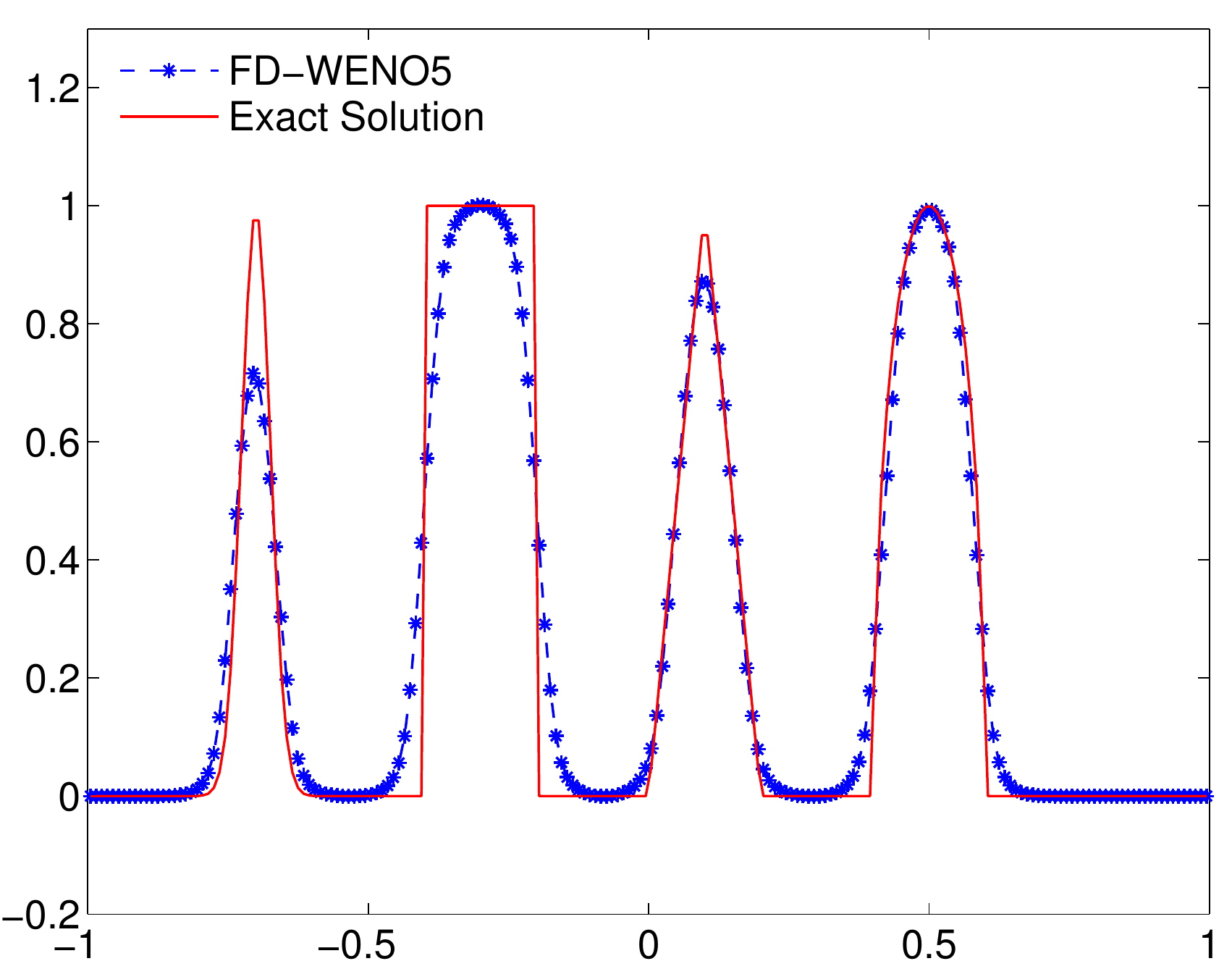} &    \includegraphics[width=6.5cm]{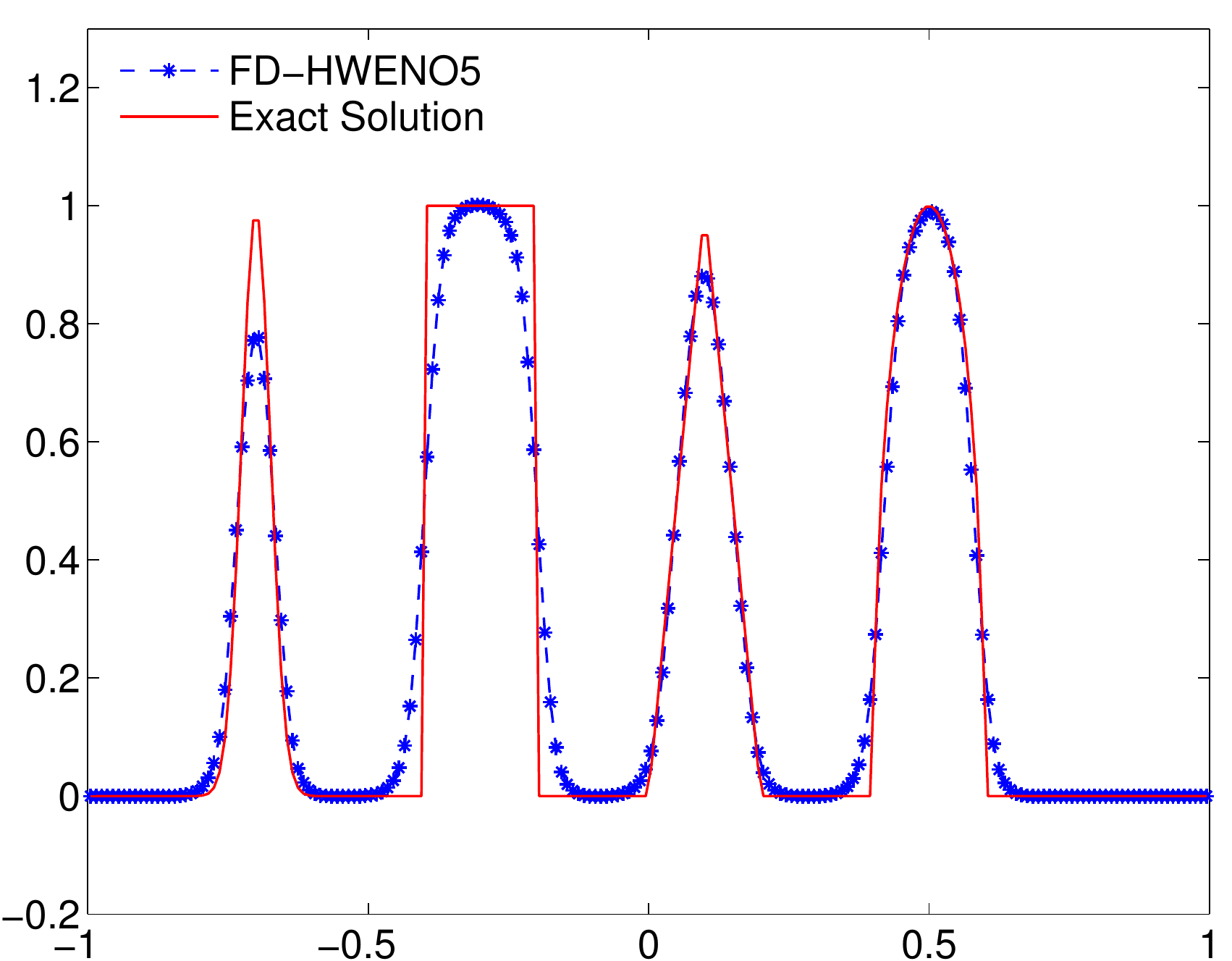} \\
  (c) Finite difference-WENO 5th                                 &    (d) Finite difference-HWENO 5th 
 \end{tabular}
\caption{\label{fig:1Dtransport}1D transport equation : {\it Plot solutions of the linear equation~\eqref{eq:test1D} with initial data~\eqref{eq:discontinuous_initial2}. $n_x=200$, CFL$=2.5$ for semi-Lagrangian methods and CFL$=0.85$ for finite difference methods. The  final time is $T_{\text{end}}=8$.}}
 \end{center}
\end{figure}

We have similar observations of both regular solution and discontinuous solution cases in Figure~\ref{fig:1Dtransport}. The  usual semi-Lagrangian method with cubic spline interpolation approximates well exponential function, but involves too much oscillation in step function. The other methods with WENO or HWENO reconstruction avoid a lot spurious oscillations. Semi-Lagrangian methods is less dissipative than finite difference method, which can be seen in step function and peak function. Moreover, the finite difference method with fifth order Hermite WENO reconstruction is less dissipative than the one with usual fifth order WENO reconstruction. 

%%%%%%%%%%%%%%%%%%%%%%%%%%%%%%%%%%%%%%%%%%%%%%%%%%%%%%%%%%%%%%%%%%%%%%%
%%                                                                   %%
%%                                                                   %%
%%%%%%%%%%%%%%%%%%%%%%%%%%%%%%%%%%%%%%%%%%%%%%%%%%%%%%%%%%%%%%%%%%%%%%%

\subsection{Simplified paraxial Vlasov-Poisson model}
\label{sec:model}

 We apply the numerical methods presented in previous sections to the following Vlasov-Poisson system satisfied by $f(t,r,v)$, where $r\in\mathbb{R}$, $v\in\mathbb{R}$~\cite{bibFS2,bibCLM}
 \begin{equation}
 \label{eq:beam}
  \left\{
   \begin{array}{l}
   \displaystyle \frac{\partial f}{\partial t} \,\,+\,\, \frac{v}{\varepsilon}\,\frac{\partial f}{\partial r} \,\,+\,\, \left(E_f - \frac{r}{\varepsilon}\right)\,\frac{\partial f}{\partial v}\,\,=\,\,0, \\[5mm]
  \displaystyle    \frac{1}{r} \frac{\partial }{\partial r}\left( r E_f \right) \,\,=\,\, \int_{\mathbb{R}} f \,dv.
   \end{array}
  \right. 
 \end{equation}
 %with the conventions $f(t,-r,-v)=f(t,r,v)$ and $E(t,-r)=-E(t,r)$.
 The electric field can be expressed explicitly as follows
 \begin{equation*}
  E_f(t,r) \,\,=\,\, \frac{1}{r}\int_0^r s\,\rho(t,s)\,ds, 
 \end{equation*}
where $\rho(t,r)=\int_{\mathbb{R}}f(t,r,v)dv$, hence we will compute $E_f$ by a simple numerical integration.

 The initial condition is chosen as a Gaussian in velocity multiplied  by a regularized step function in $r$:
 \begin{equation}
  f_0(r,v)=\frac{4}{\sqrt{2\pi\alpha}}\chi(r)\exp(-\frac{v^2}{2\alpha}), 
 \end{equation}
 with $\chi(r)=\frac{1}{2}\text{erf}(\frac{r+1.2}{0.3})-\frac{1}{2}\text{erf}(\frac{r-1.2}{0.3})$ and  $\alpha=0.2$.  The Vlasov-Poisson system~\eqref{eq:beam} conserves mass
 \begin{equation*}
  \frac{d   }{dt} \int_{\mathbb{R}^2} f(t,r,v) dr dv = 0,
 \end{equation*}
and also $L^p$ norm for $1\leq p < \infty$
 \begin{equation*}
  \frac{d   }{dt} ||f(t,r,v)||_{L^p(\mathbb{R}^2)} = 0.
 \end{equation*}
Therefore, the evolution in time of these quantities will be observed for various numerical schemes. We will also investigate the time evolution of the kinetic energy of the Vlasov-Poisson system~\eqref{eq:beam} :
\begin{equation}
 \mathcal{E}(t)\,=\,\int_{\mathbb{R}^2}\frac{v^2}{2}\,f \,dr\, dv.
\end{equation}
A reference solution of kinetic energy is computed using a fifth order finite difference WENO method with very fine mesh ($n_x=1025$, $\Delta t=1/1600$).

In the following we take the parameter $\varepsilon=0.7$ and the computational  domain is $(r,v)\in\Omega=[-4,4]^2$.

%%------------------------------------------------------------------
Concerning the numerical resolution using  semi-Lagrangian methods,  we notice that we deliberately choose not to apply a time splitting in order to use this method in a general context.  The characteristic curves corresponding to the Vlasov equation~\eqref{eq:beam} cannot be solved explicitly. Then we apply a second order leap-frog scheme already developed in~\cite{bibSR}. Finally,  we interpolate the distribution function $f(r^n,v^n)$ by a tensor product for cubic spline or  by a dimension by dimension Hermite WENO reconstruction.

\begin{figure}
\begin{center}
 \begin{tabular}{cc}
\includegraphics[width=6.5cm]{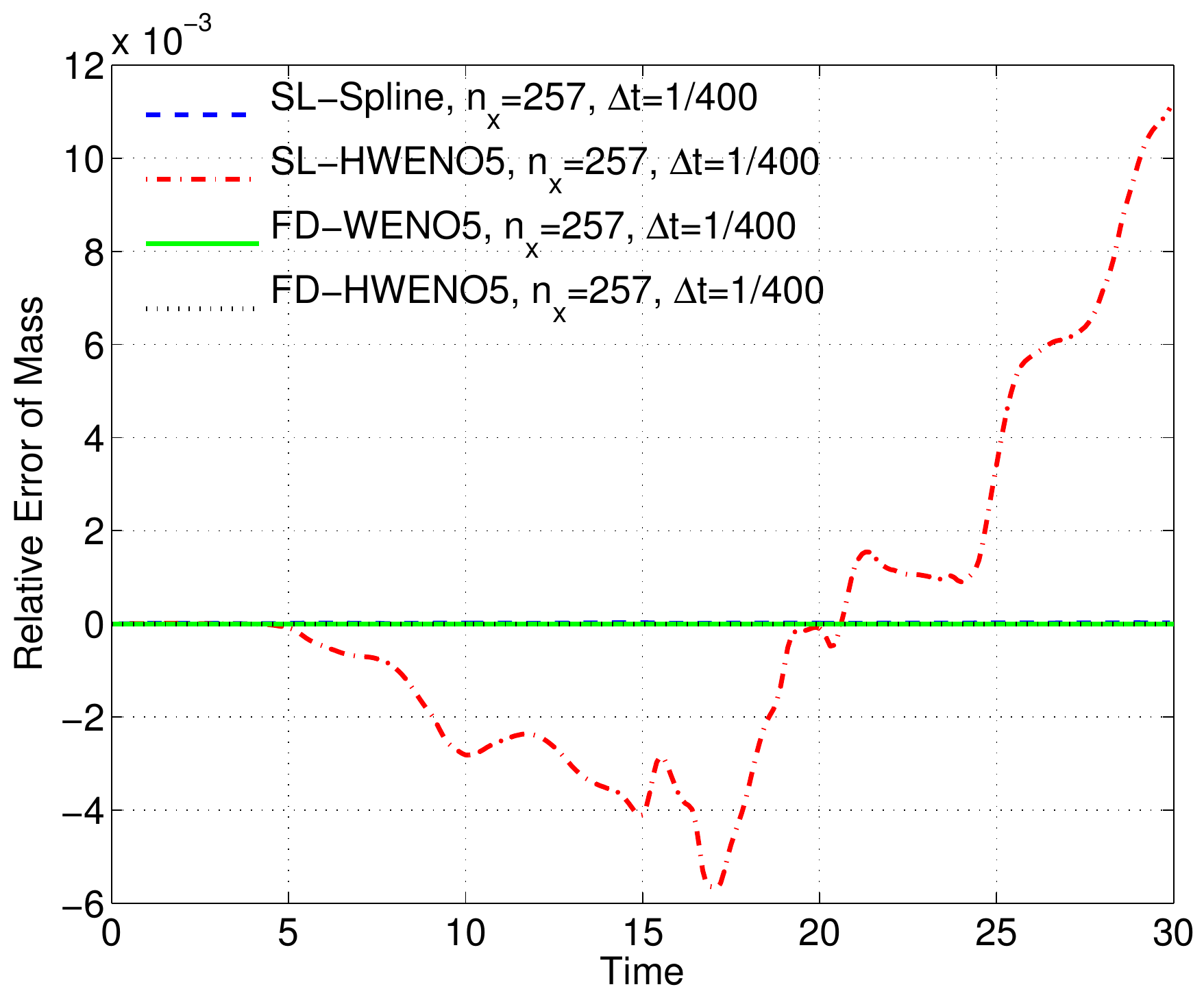} &    
\includegraphics[width=6.5cm]{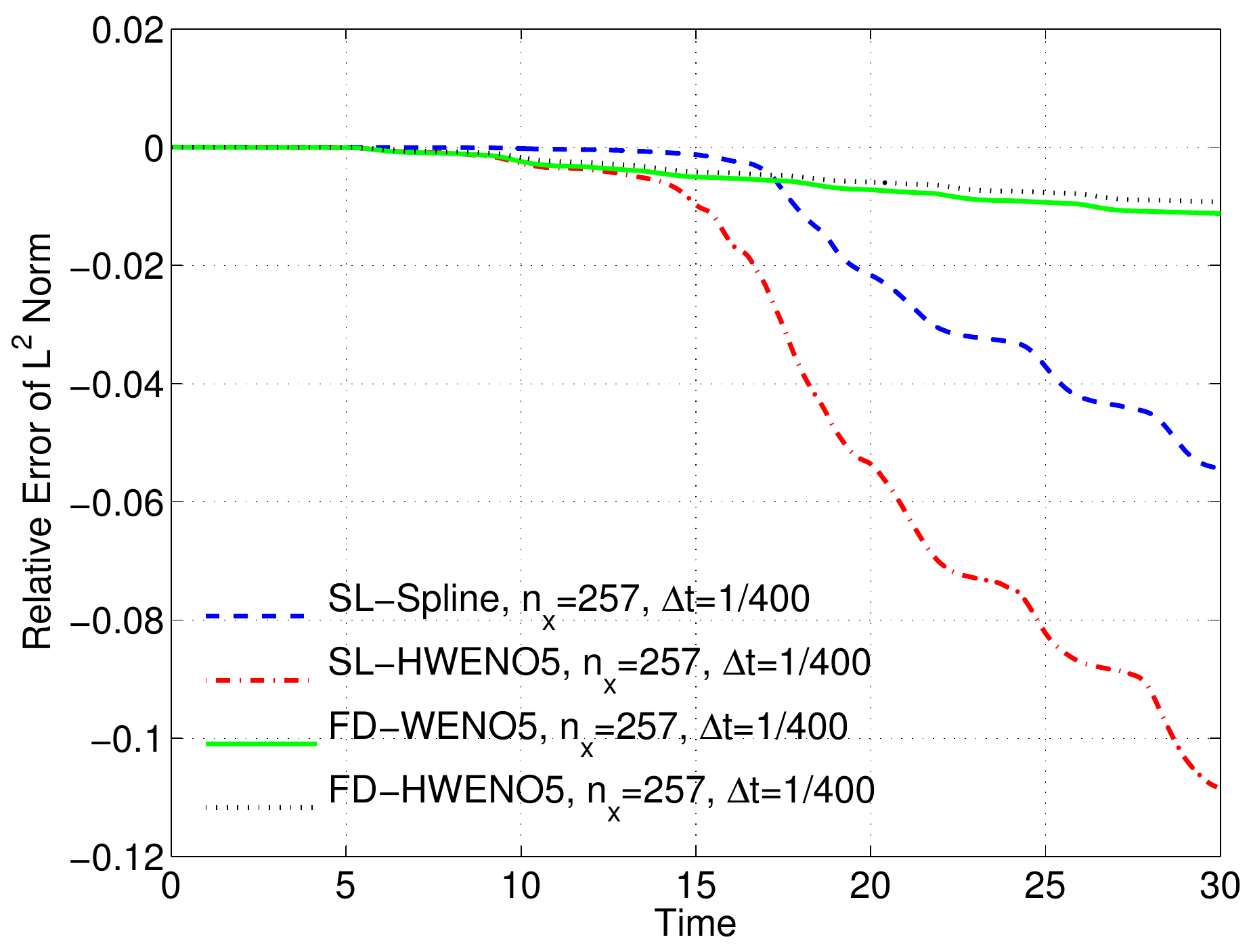}
\\
 (a) &    (b)              
\\
\includegraphics[width=6.5cm]{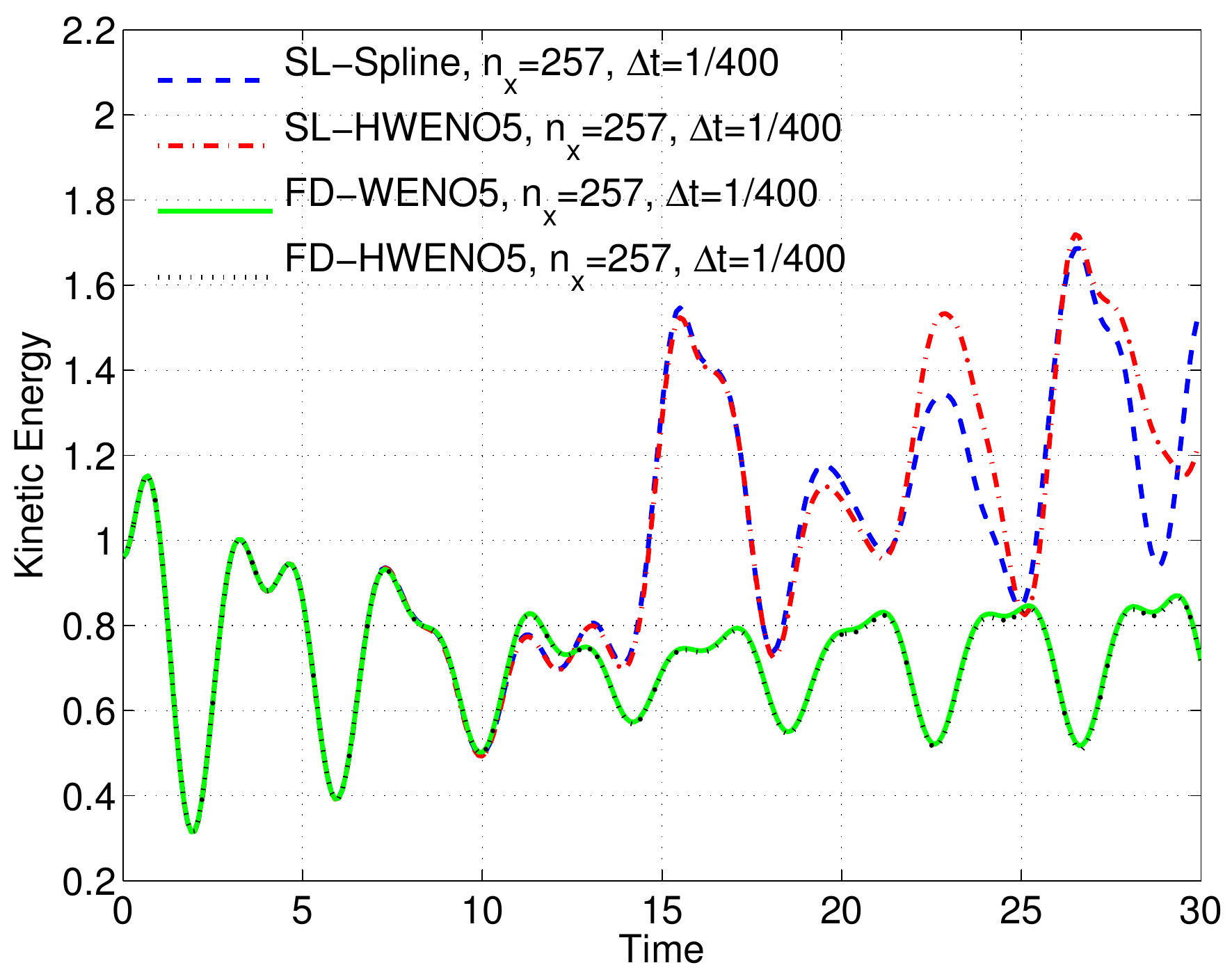}&
\includegraphics[width=6.5cm]{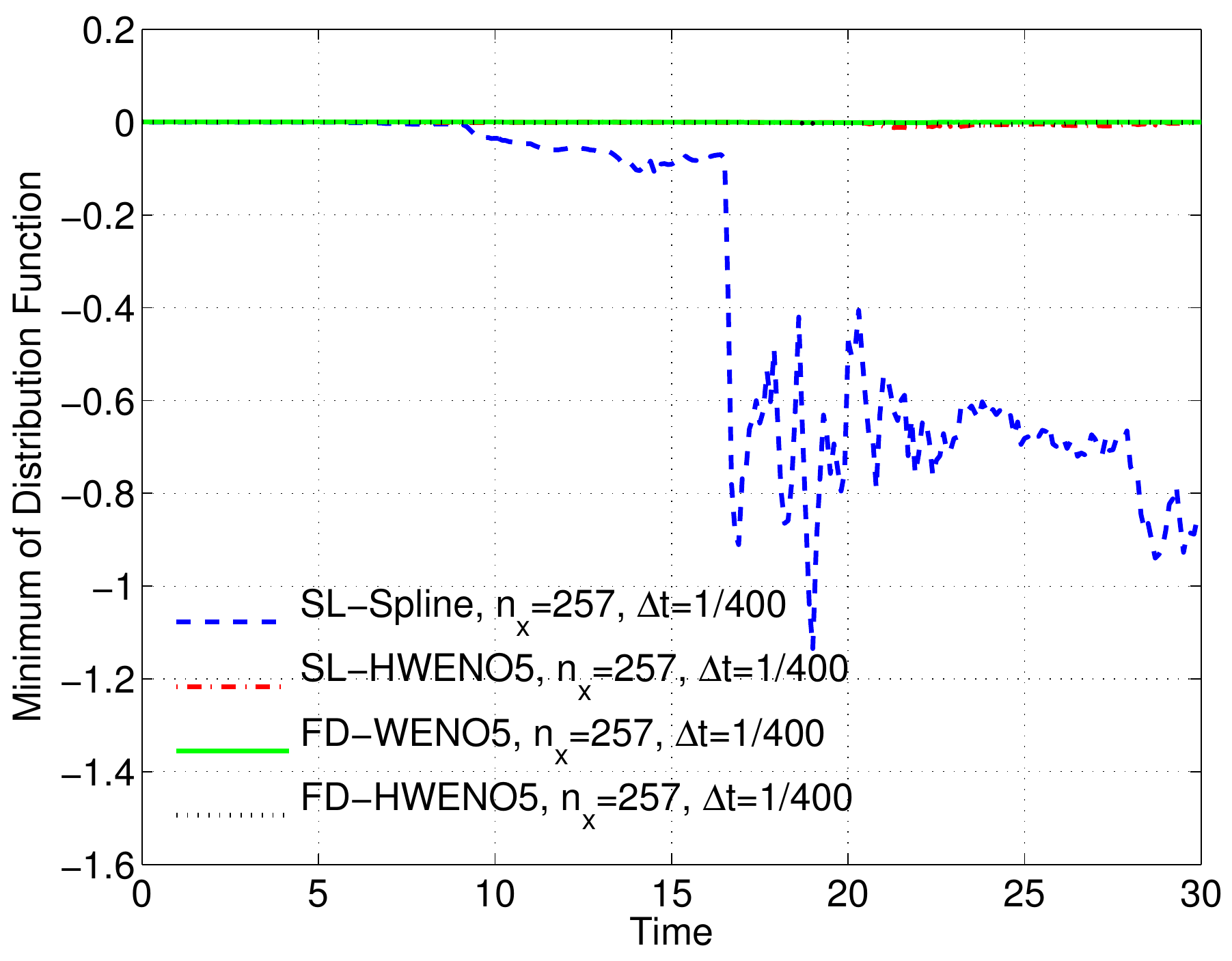}
\\
 (c) &    (d) 
 \end{tabular}
\caption{\label{fig:Beam_spline_hweno5}Simplified paraxial Vlasov-Poisson model : {\it Comparison between semi-Lagrangian with cubic spline, Hermite WENO 5th and finite difference with Hermite WENO 5th methods for Beam test.}}
 \end{center}
\end{figure}

In Figure~\ref{fig:Beam_spline_hweno5} we compare the evolution of invariants (mass, $L^2$ norm, minimum f the density) and the kinetic energy obtained from   semi-Lagrangian and finite difference methods  to the reference solution. We denote  the linear phase for time interval $t\leq 10$, and   the nonlinear phase for time interval $t>10$ where small filaments are generated.
 
As we can see from Figure~\ref{fig:Beam_spline_hweno5}, the kinetic energy obtained with the semi-Lagrangian method with cubic spline and Hermite WENO5 reconstruction  is relatively close to the reference solution in the linear phase, but it diverges from the  reference one in the nonlinear phase.  Even with a fines mesh, we cannot improve the numerical results for large time. 

\begin{figure}
\begin{center}
 \begin{tabular}{cccc}
\includegraphics[width=3cm]{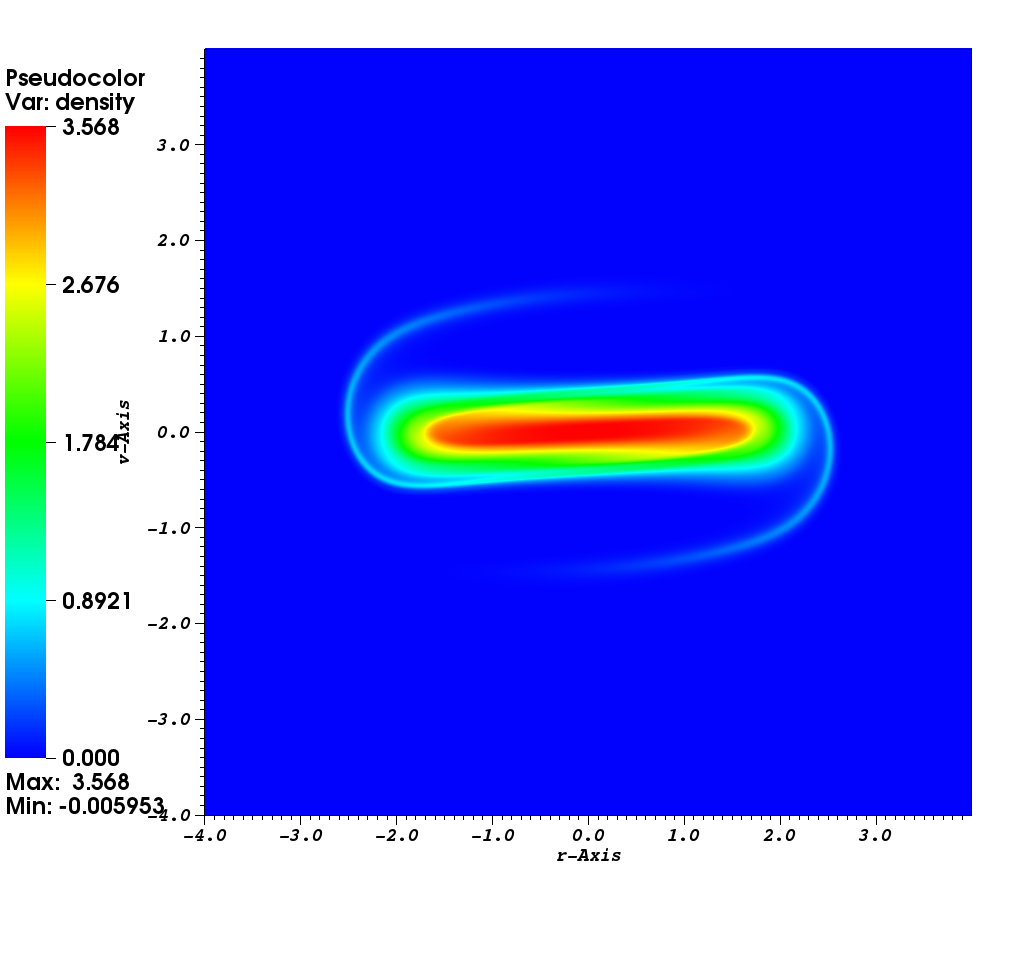} & 
\includegraphics[width=3cm]{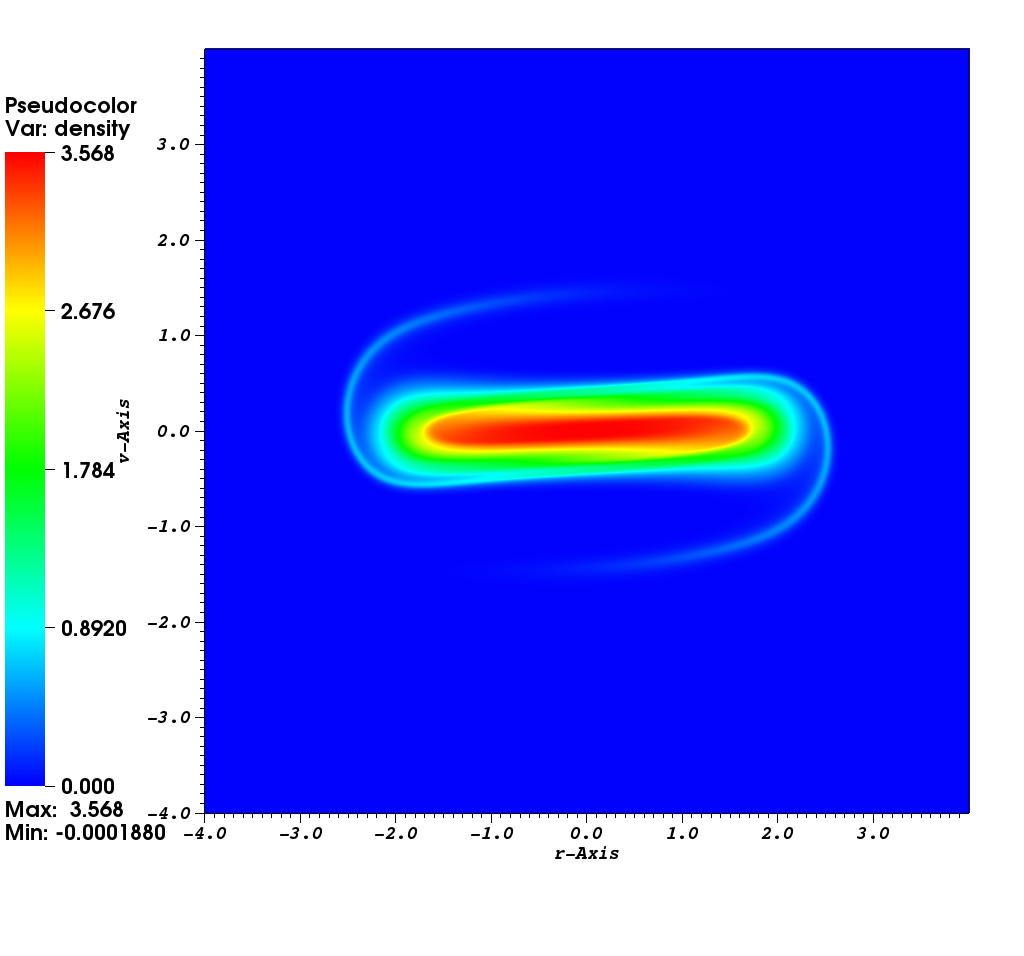} &
\includegraphics[width=3cm]{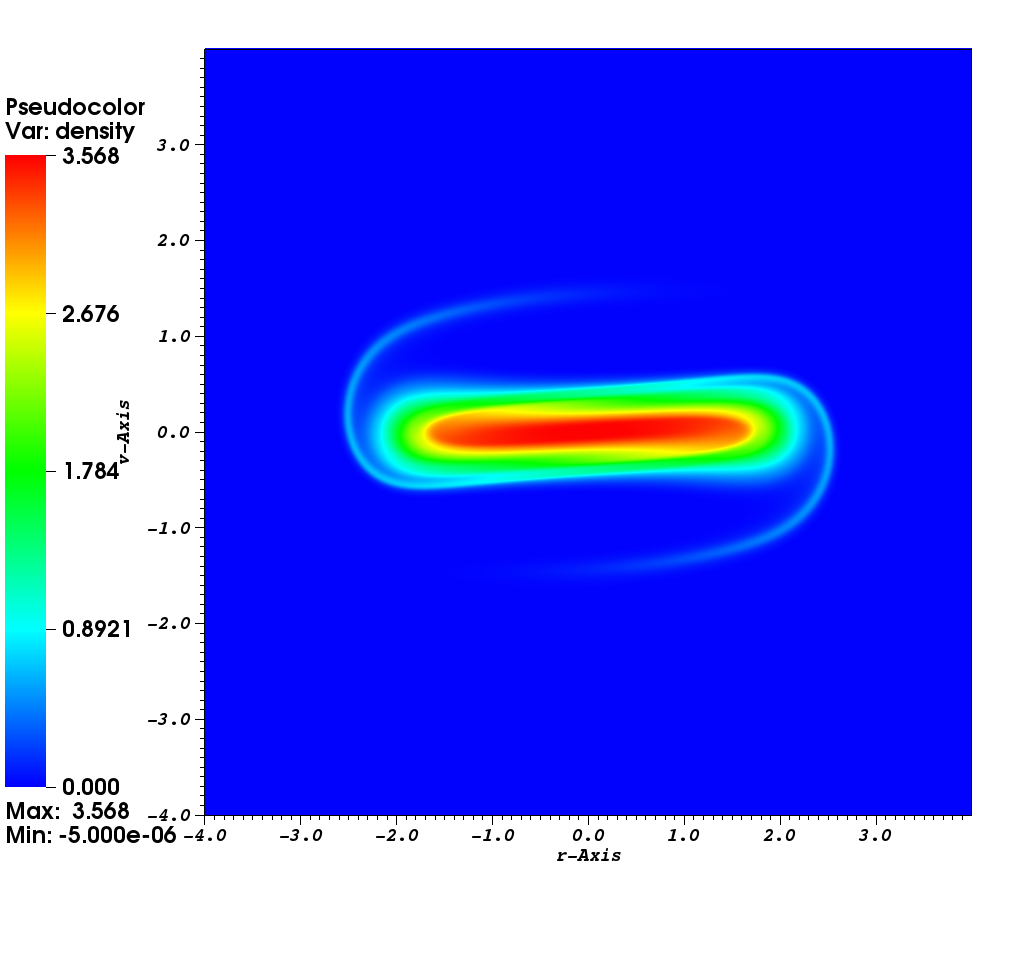} &    
\includegraphics[width=3cm]{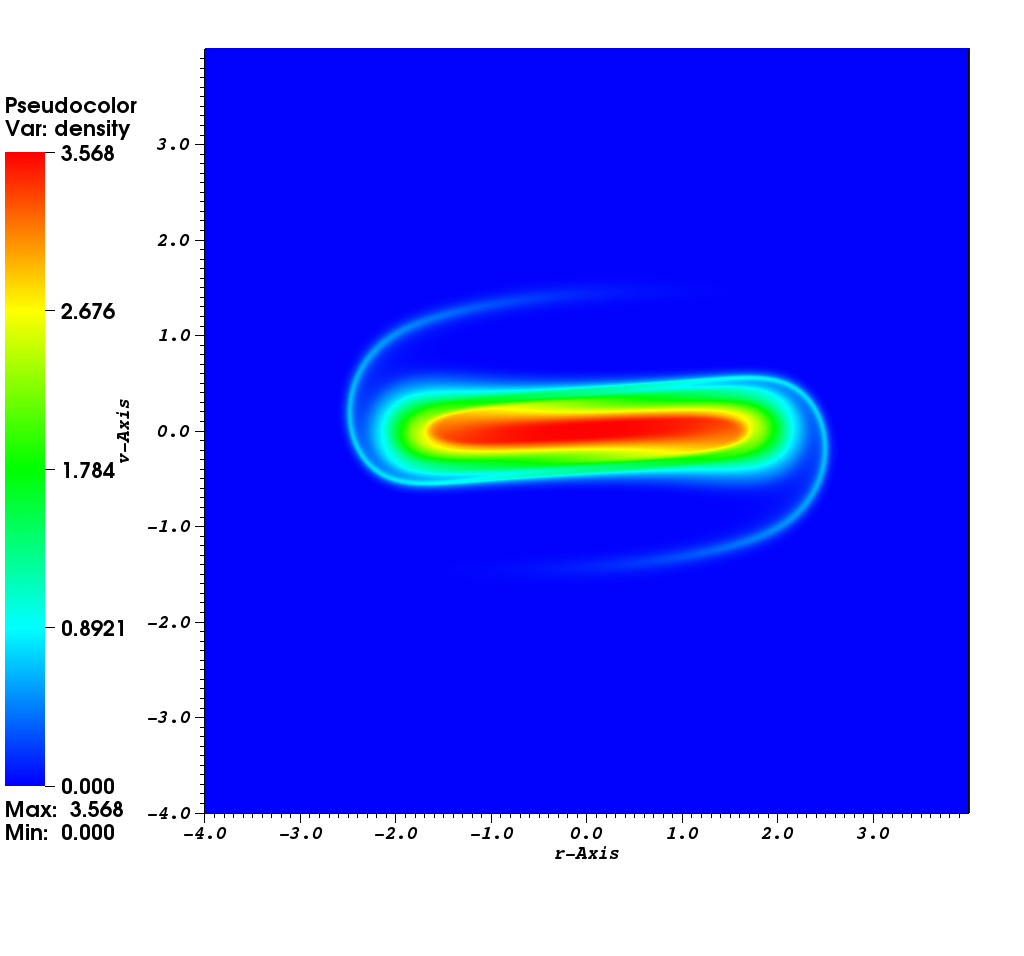} 
\\
 % (a) $t=10$ &    (b)   $t=10$ &    (c)   $t=10$  & (d)   $t=10$  
%\\
\includegraphics[width=3cm]{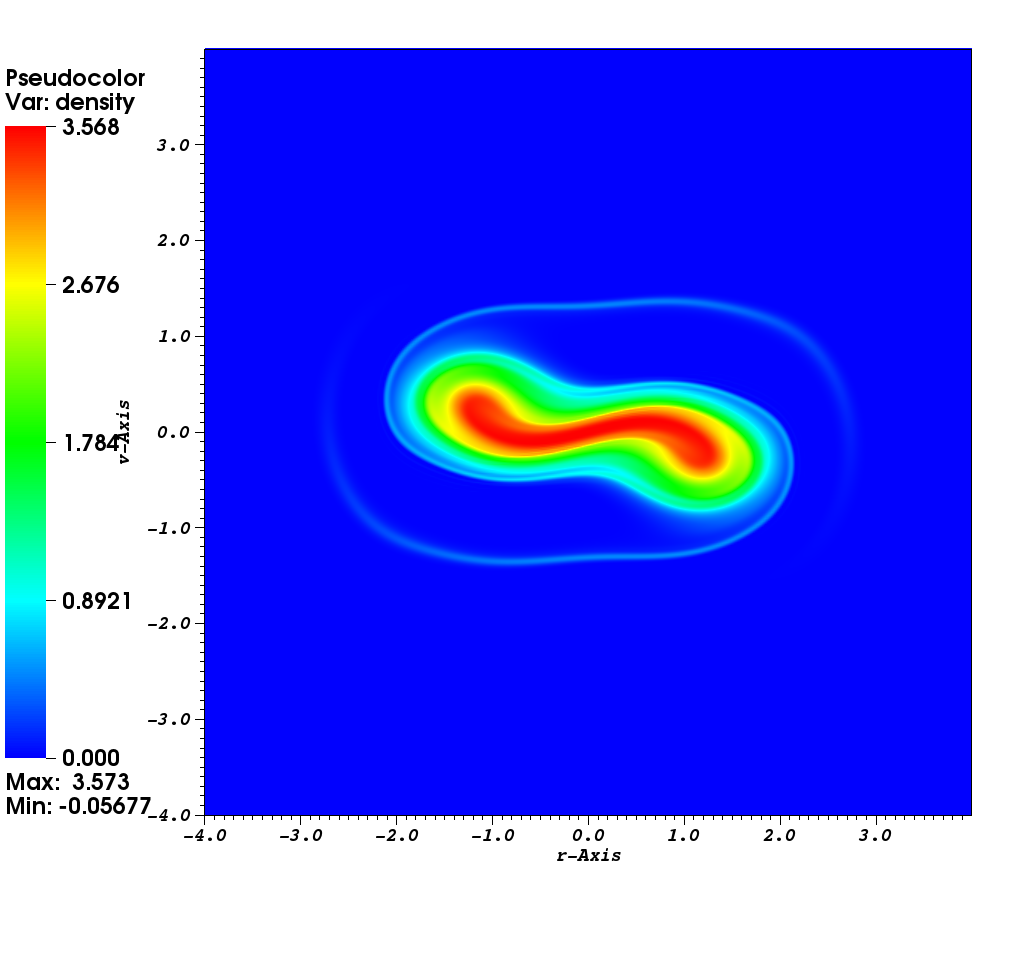} & 
\includegraphics[width=3cm]{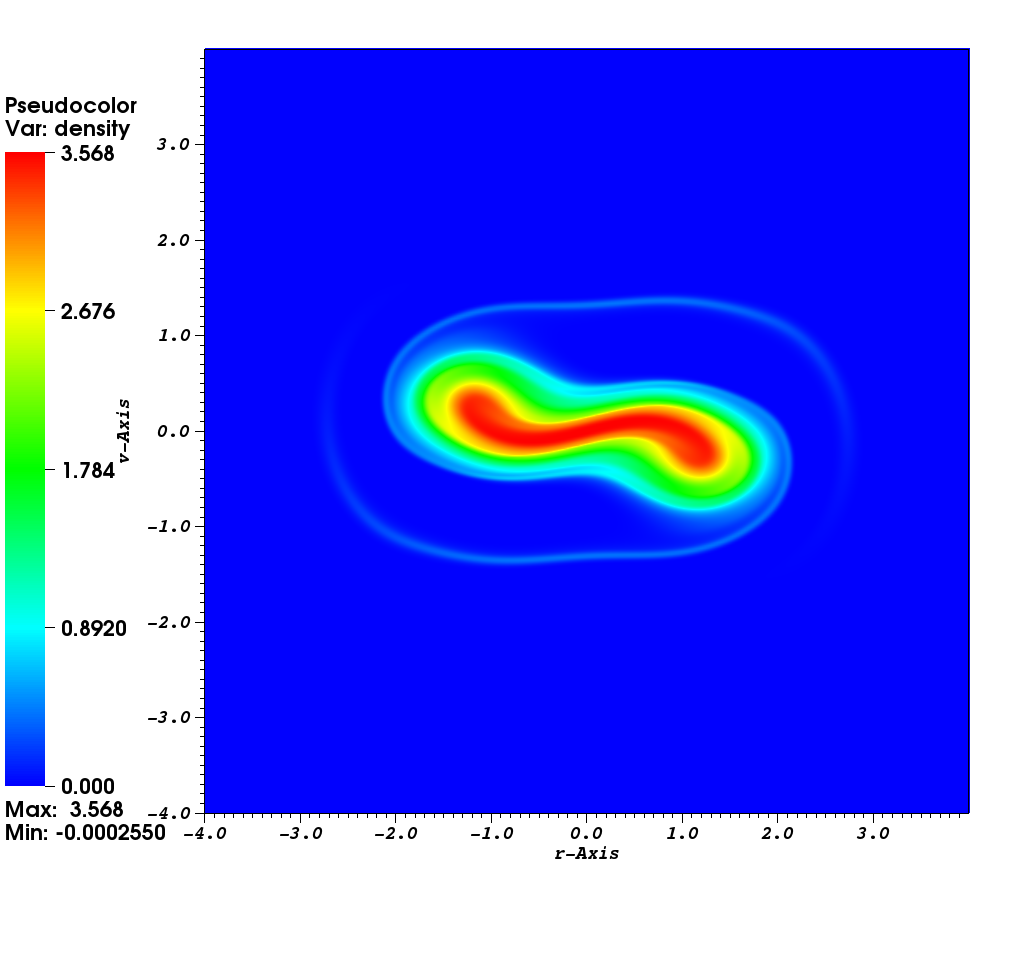} &  
\includegraphics[width=3cm]{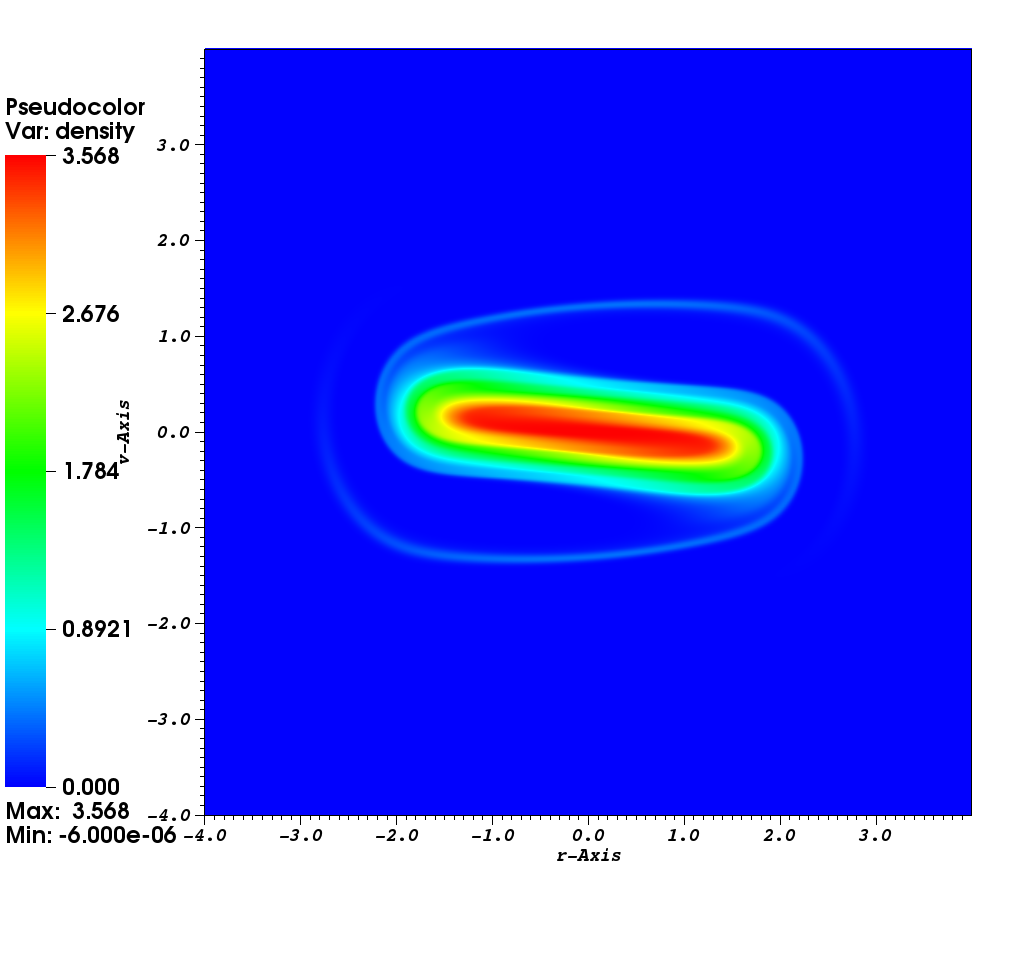} &
\includegraphics[width=3cm]{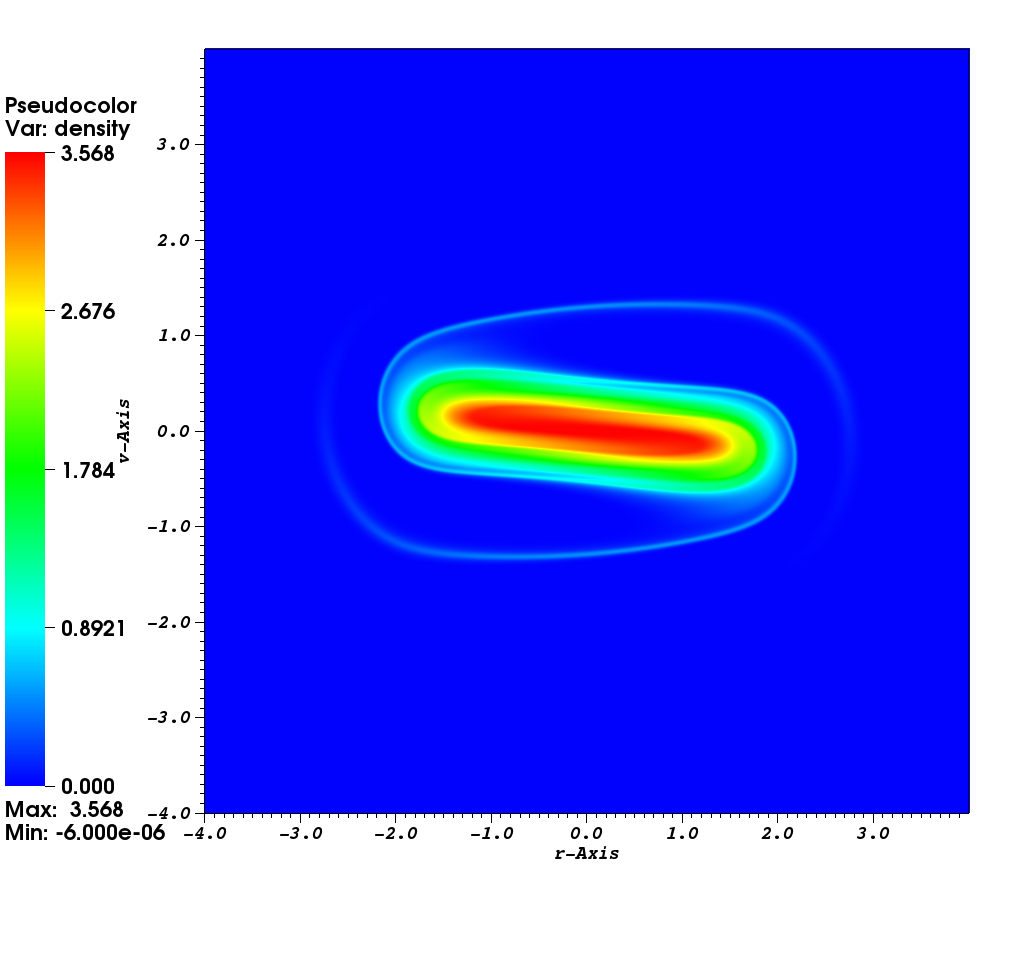} \\
%(e)$t=15$ & (f)  $t=15$ & (g)  $t=15$ & (h)  $t=15$ 
%\\
\includegraphics[width=3cm]{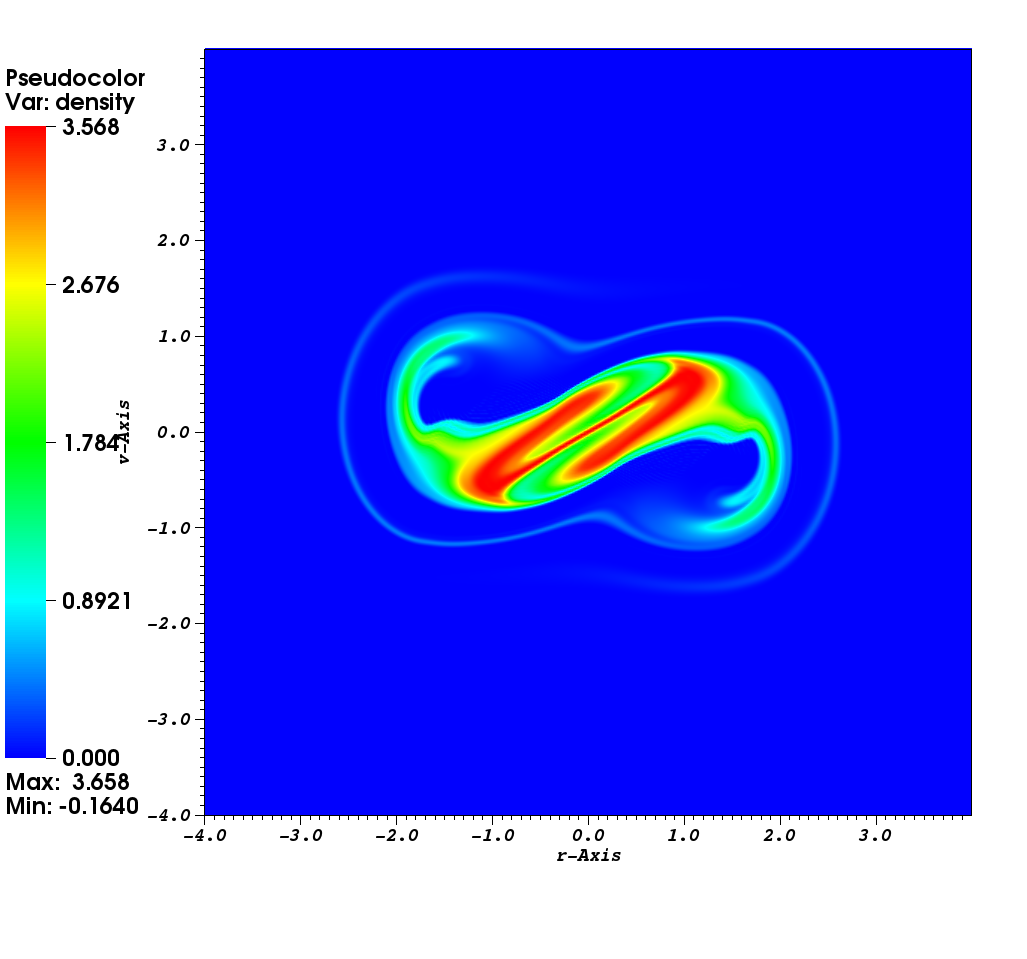}  & 
\includegraphics[width=3cm]{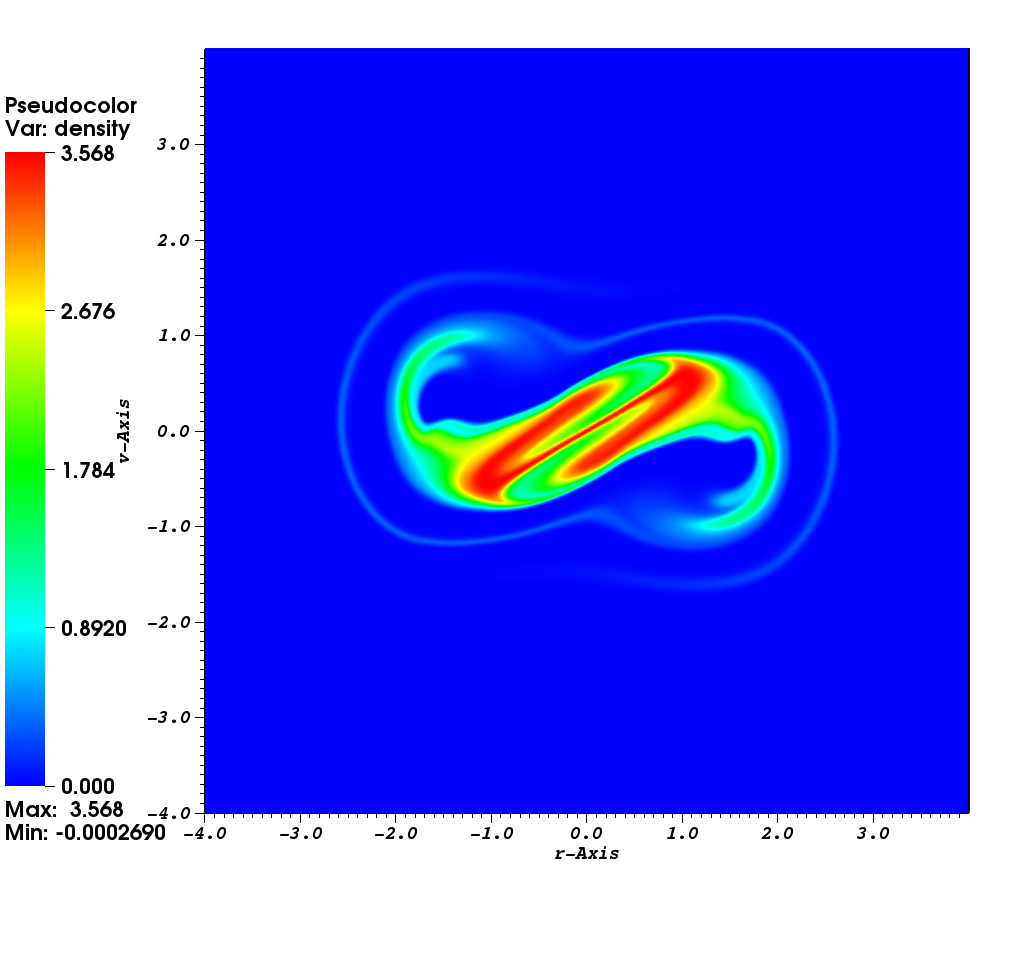} &    
\includegraphics[width=3cm]{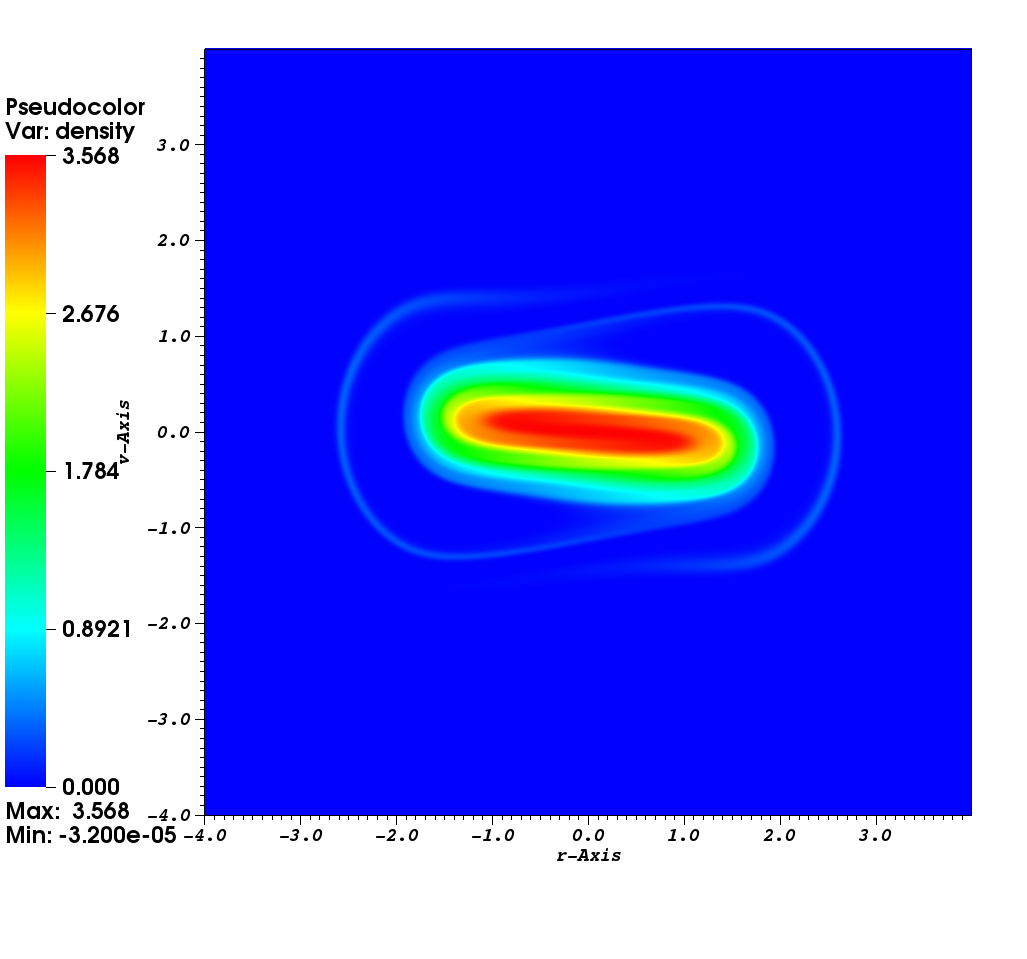} &
\includegraphics[width=3cm]{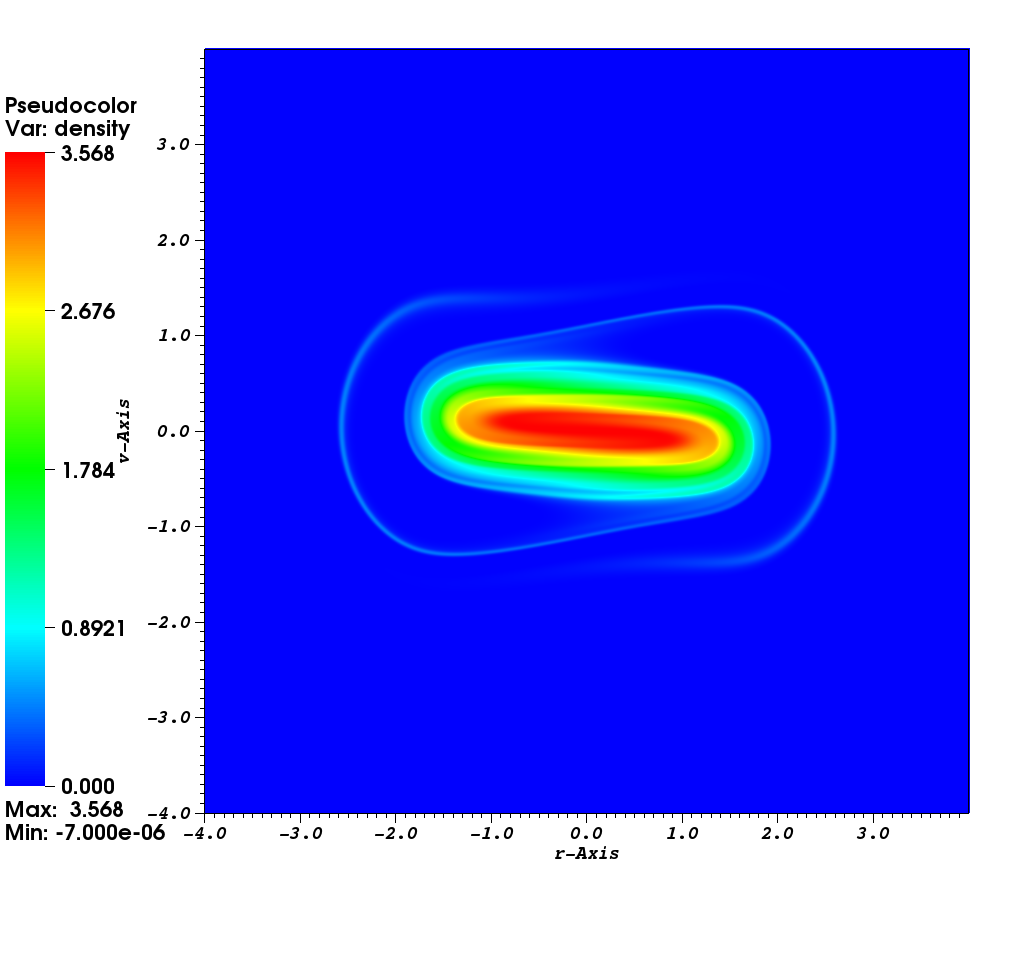} \\
%  (i) $t=20$&    (j)   $t=20$ &    (k)   $t=20$ & (l)$t=20$ \\
 (a)  &    (b)    &    (c)    & (d)  

 \end{tabular}
\caption{\label{fig:Beam_SL_distribution}Simplified paraxial Vlasov-Poisson model : {\it Distribution function for Beam test: (a) semi-Lagrangian with cubic spline; (b)  semi-Lagrangian with Hermite WENO5; (c) finite difference with Hermite WENO5; (d) reference solution at time $t=10$, $15$ and $20$. Mesh size is  $n_x=513, \Delta t=1/800$.}}
 \end{center}
\end{figure}

Finally, we compare the distribution function $f$ obtained from  semi-Lagrangian with cubic spline and Hermite WENO5 methods and finite difference with Hermite WENO5 reconstructions  with a reference solution computed with a refined mesh ($\Delta t=1/1600$ and $n_x=1025$) in Figure~\ref{fig:Beam_SL_distribution}.
At time $t=10$, the distribution function $f$ of  semi-Lagrangian methods is very close to the  reference solution, where two small filaments appear (see top of Figure~\ref{fig:Beam_SL_distribution}). Then during the nonlinear phase, a large number of filaments are generated due to the non-linearity of the Vlasov-Poisson system and the distribution functions $f$  obtained with semi-Lagrangian and finite difference methods start to differ strongly at $t\geq 15$, which also correspond to the divergence of the kinetic energy of Figure~\ref{fig:Beam_spline_hweno5}. At time  $t\geq 20$, the semi-Lagrangian methods generates a completely unstable beam who is not consistent with the results obtained for the reference solution.

These numerical simulations illustrate perfectly that  semi-Lagrangian methods without splitting  work well during the linear phase even with very large time step, but  they do not seem  very robust during the nonlinear phase when micro-structures appear (filamentation).

%%%%%%%%%%%%%%%%%%%%%%%%%%%%%%%%%%%%%%%%%%%%%%%%%%%%%%%%%%%%%%%%%%%%%%%
%%                                                                   %%
%%                                                                   %%
%%%%%%%%%%%%%%%%%%%%%%%%%%%%%%%%%%%%%%%%%%%%%%%%%%%%%%%%%%%%%%%%%%%%%%%

\subsection{Guiding center model}
\label{sec:Num}

%In the previous section, we have seen that semi-Lagrangian method works well in linear phase and it is more efficient in computational cost. As contrast, finite difference methods should respect CFL condition, but they are more stable in non-linear phase. We thus combine these two method for plasma physical simulation.

We finally consider the  guiding center model~\cite{bibCLM}, which has been derived to describe highly magnetized plasma in the transverse plane
\begin{equation}
\left\{
\begin{array}{l}
 \frac{\partial \rho}{\partial t}+\mathbf{U}\cdot\nabla\rho=0,\\[3mm]
 -\Delta\phi=\rho.
\end{array}
\right.
   \label{eq:guiding_center}
  \end{equation}
where the velocity $\mathbf{U} = (-\partial_y\phi, \partial_x\phi)  $. Here we consider the model in a disk domain
\begin{equation*}
 D\,\,=\,\,\{(x,y)\in\mathbb{R}^2 : \sqrt{x^2+y^2}\leq R\}
\end{equation*}
and assume that the electric potential is vanishing at the boundary
\begin{equation}
 \phi(x,y)=0,\quad(x,y)\in\partial D.
 \label{eq:GC_BC}
\end{equation}
Then if we ignore the effect of boundary conditions, the  guiding center model verifies the following properties :
\begin{enumerate}
%%---------------------------------------------------------
 \item Positivity of density $\rho$
 \begin{equation*}
  0\leq\rho(t,x,y).  
 \end{equation*}
  %%---------------------------------------------------
 \item Mass conservation
 \begin{equation*}
  \frac{d}{dt}\left(\int_{D}\rho dx dy\right)=0.
 \end{equation*}
 %%---------------------------------------------------
 \item $L^p$ norm conservation, for $1\leq p\leq\infty$
 \begin{equation*}
  \frac{d}{dt}||\rho||_{L^p(D)}=0.
 \end{equation*}
%%-----------------------------------------------------
\item Energy conservation
\begin{equation*}
 \frac{d}{dt}\left(\int_{D}|\nabla\phi|^2dx dy\right)=0.
\end{equation*}
\end{enumerate}

To solve the system~\eqref{eq:guiding_center}, we use a scheme based on Cartesian mesh and apply an Inverse Lax-Wendroff procedure to treat boundary conditions on kinetic equations \cite{bibFY1,bibFY2}. Actually, this system has been already solved in polar coordinates in \cite{bibHMMP}. But, the change of coordinate usually produce an artificial  singularity at the origin, which needs a particular treatment.  At contrast, with the Inverse Lax-Wendroff technique on Cartesian mesh, we do not have any singularity and it is not related to the numerical scheme since boundary effects and numerical schemes are treated independently. Furthermore,  it is easy to adapt to other geometries~\cite{bibFY1, bibFY2}.

In this section, we only focus on discretization of boundary condition~\eqref{eq:GC_BC} of Poisson equation.  The one for transport equation is trivial, since a homogeneous Dirichlet boundary condition will be used.

%%%%%%%%%%%%%%%%%%%%%%%%%%%%%%%%%%%%%%%%%%%%%%%%%%%%%%%%%%%%%%%%%%%%%%%
%%                                                                   %%
%%                                                                   %%
%%%%%%%%%%%%%%%%%%%%%%%%%%%%%%%%%%%%%%%%%%%%%%%%%%%%%%%%%%%%%%%%%%%%%%%
\subsubsection{Discretization of  Poisson equation}
%%%%%%%%%%%%%%%%%%%%%%%%%%%%%%%%%%%%%%%%%%%%%%%%%%%%%%%%%%%%%%%%%%%%%%%
\label{sec:discratisation:parabolique}

\begin{figure}[h]
  \begin{center} 
    \includegraphics[width=10cm]{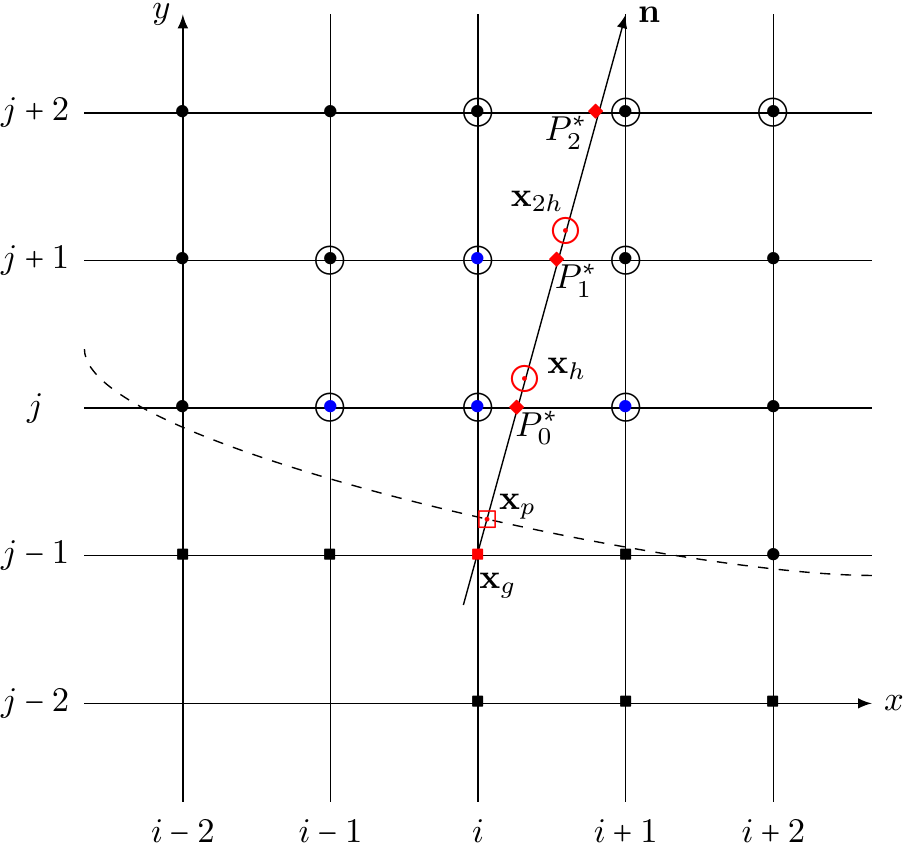}
\caption{\label{fig:2Ddomain}Spatially two-dimensional Cartesian mesh. $\bullet$ is interior point, $\filledsquare$ is ghost point, $\boxdot$ is the point at the boundary, $\largecircle$ is the point for extrapolation, the dashed line is the boundary.}
  \end{center}
\end{figure}

A classical five points finite difference approximation is used to discretize  the Poisson equation. 
However, to discretize the Laplacian operator $\Delta\phi$ near the physical boundary, we notice that some points of the usual five points finite difference formula can be located outside of interior domain. 
For instance,   Figure~\ref{fig:2Ddomain} illustrates the discretization stencil for  $\Delta\phi$ at the point $(x_i,y_j)$.
The point $\mathbf{x}_g=(x_i,y_{j-1})$ is located outside of interior domain. Let us denote the approximation of $\phi$ at  the point $\mathbf{x}_g$ by $\phi_{i,j-1}$.
Thus $\phi_{i,j-1}$ should be extrapolated from the interior domain.

We extrapolate $\phi_{i,j-1}$ on the normal direction $\mathbf{n}$ 
\begin{equation}
 \phi_{i,j-1}=\tilde{w}_p {\phi}(\mathbf{x}_p) + \tilde{w}_h {\phi}(\mathbf{x}_h) + \tilde{w}_{2h} {\phi}(\mathbf{x}_{2h}) ,
 \label{eq:normal_extrapolation}
\end{equation}
where $\mathbf{x}_p$ is the cross point of the normal $\mathbf{n}$ and the physical boundary $\Omega_{\mathbf{x}_\bot}$.
The points $\mathbf{x}_h$ and $\mathbf{x}_{2h}$ are  equal spacing on the normal $\mathbf{n}$, {\it i.e.} $h=|\mathbf{x}_p-\mathbf{x}_h|=|\mathbf{x}_h-\mathbf{x}_{2h}|$,
with $h=\min(\Delta x,\Delta y)$, $\Delta x$, $\Delta y$ are the space steps in the directions $x$ and $y$ respectively.
Moreover, $\tilde{w}_p$, $\tilde{w}_h$, $\tilde{w}_{2h}$ are the extrapolation weights depending on the position of $\mathbf{x}_g$, $\mathbf{x}_p$, $\mathbf{x}_h$ and $\mathbf{x}_{2h}$.
In~\eqref{eq:normal_extrapolation}, ${\phi}(\mathbf{x}_p)$ is given by the  boundary condition~\eqref{eq:GC_BC},
whereas ${\phi}(\mathbf{x}_h)$, ${\phi}(\mathbf{x}_{2h})$ should be determined by interpolation.

For this, we first construct an interpolation stencil $\mathcal{E}$, composed of grid points of $\Omega$. 
For instance, in Figure~\ref{fig:2Ddomain},  the inward normal $\mathbf{n}$ intersects  the grid lines $y=y_{j}$, $y_{j+1}$, $y_{j+2}$ at points $P^*_0$, $P^*_1$, $P^*_2$. 
Then we choose the three nearest points  of the cross point $P^*_l,\,l=0,1,2$, in each line, {\it i.e.} marked by a large circle. 
From these nine points, we construct a  Lagrange polynomial $q_2(\mathbf{x})\in\mathbb{Q}_2(\mathbb{R}^2)$ and evaluate the polynomial $q_2(\mathbf{x})$ at $\mathbf{x}_h$ and $\mathbf{x}_{2h}$, {\it i.e.}
\begin{eqnarray*}
 {\phi}(\mathbf{x}_h)&=&\sum_{\ell=0}^8w_{h,\ell}{\phi}(\mathbf{x}_{\ell}),\\[3mm]
 {\phi}(\mathbf{x}_{2h})&=&\sum_{\ell=0}^8w_{2h,\ell}{\phi}(\mathbf{x}_{\ell}),
\end{eqnarray*}
with $\mathbf{x}_{\ell}\in \mathcal{E}$. 
Hence,  we get $\phi_{i,j-1}$ which is approximated from the interior domain.

However, in some cases, we can not find a stencil of nine interior points.  For instance, when the interior domain has small acute angle sharp, the normal $\mathbf{n}$ can not have three cross points $P^*_l,\,l=0,1,2$ in interior domain, or we can not have three nearest points  of the cross point $P^*_l,\,l=0,1,2$, in each line. 
In this case, we alternatively use a first degree polynomial $q_1(\mathbf{x})$ with a four points  stencil or even a zero degree polynomial $q_0(\mathbf{x})$ with  an one point  stencil.  We can similarly construct the four points stencil or the one point stencil as above.

%%---------------------------------------------------------------------------
\subsubsection{Numerical simulation of the diocotron instability}

We now consider the diocotron instability for an annular electron layer. This  plasma instability is created by two sheets of charge slipping past each other and is the analog of the Kelvin-Helmholtz instability in fluid mechanics. The initial data is given by
\begin{equation*}
 \rho_0(\mathbf{x}_\bot)=
 \left\{
 \begin{array}{ll}
  (1+\varepsilon\cos(\ell\theta))\exp{(-4(r-6.5)^2)},&\text{if  } r^-\leq\sqrt{x^2+y^2}\leq r^+,\\[3mm]
  0,&\text{otherwise},
 \end{array}
 \right.
\end{equation*}
where $\varepsilon$ is a small parameter, $\theta=\text{atan2}(y,x)$.
In the following tests, we take $\varepsilon=0.001$, $r^-=5$, $r^+=8$, $\ell=7$.
%%%%%%%%%%%%%%%%%%%%%%%%%%%%%%%%%%%%%%%%%%%%
%
%%%%%%%%%%%%%%%%%%%%%%%%%%%%%%%%%%%%%%%%%%%%

We have seen in the previous section that semi-Lagrangian methods may be not very appropriate. Indeed, the semi-Lagrangian method has some limitations during the  nonlinear phase,  when small filaments appear since a small time step must be used and the method is no more conservative. Therefore, we propose to apply a mixed method based on the Hermite interpolation with a WENO reconstruction: we use the semi-Lagrangian method for the linear phase with large time step; then we  apply the conservative finite difference scheme for nonlinear phase with small time step respecting CFL condition. The criterion to pass from semi-Lagrangian to finite difference methods is as follows
\begin{equation}
 \left|\int_{\RR^2}\left[\rho_h(t_n)-\rho_h(t_{n-1})\right]d\xx \right|\,>\, h^3,
\end{equation}
where $h$ is the smallest space step. 

%We remark that $h^3$ is chosen instead of $h^5$.  It is because the parabolic assumption method~\cite{bibGV} is a second order method, thus we can not achieve fifth order convergence even with a regular initial data.

\begin{figure}
\begin{center}
 \begin{tabular}{cc}
\includegraphics[width=6.5cm]{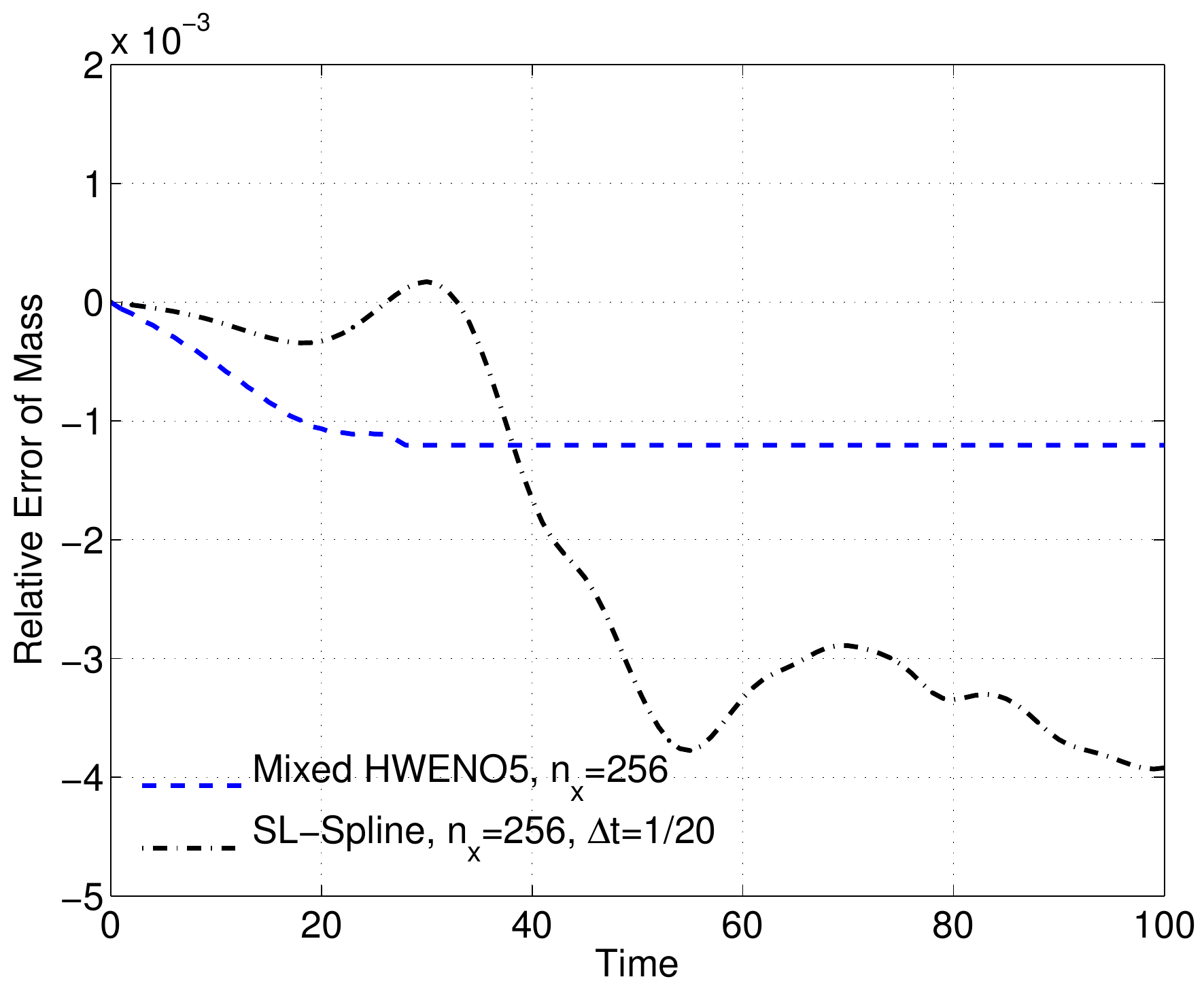}   &    
\includegraphics[width=6.5cm]{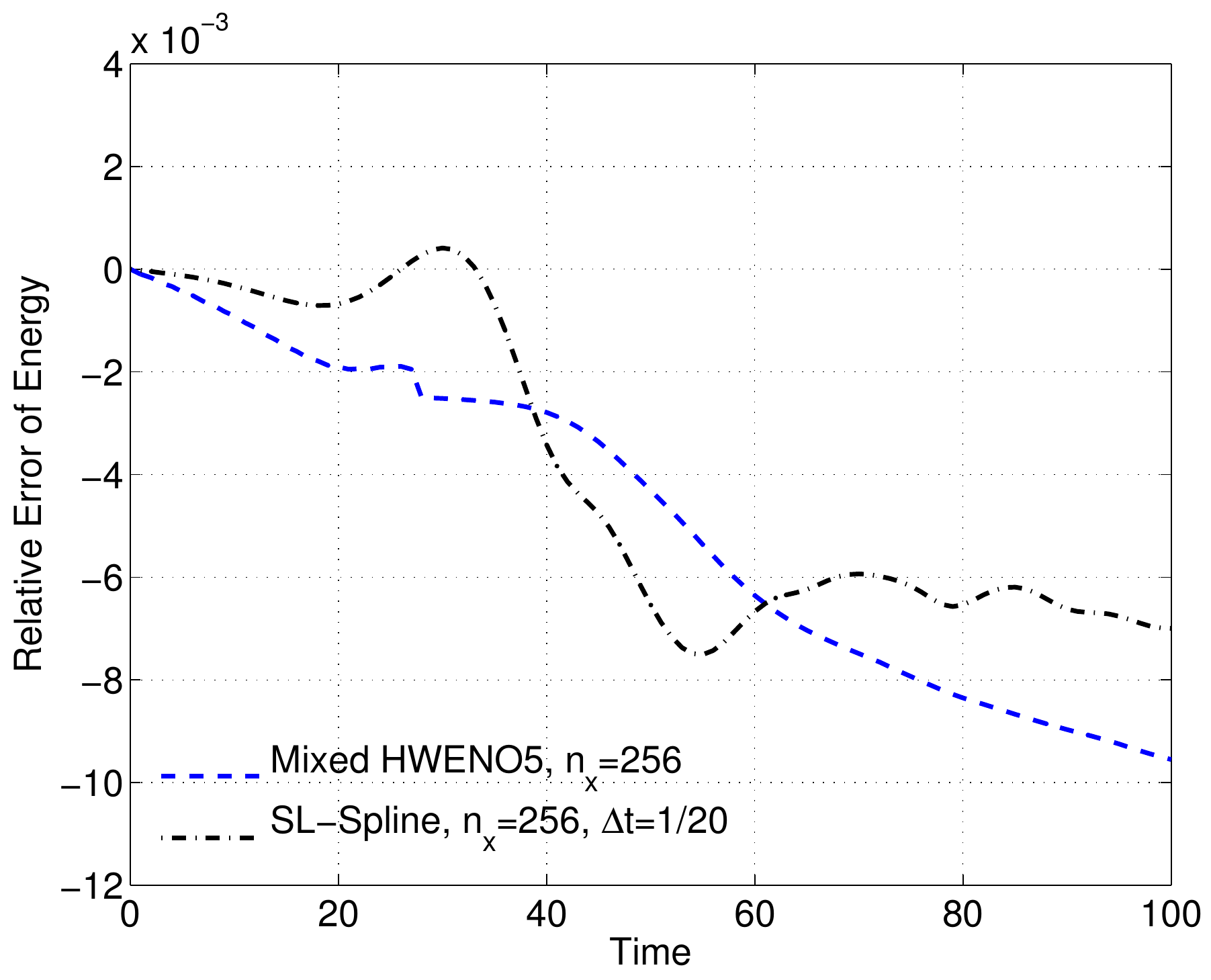}
\\
  (a)       Relative Error of Mass                              &   
 (b)  Relative Error of Energy  \\

  \includegraphics[width=6.5cm]{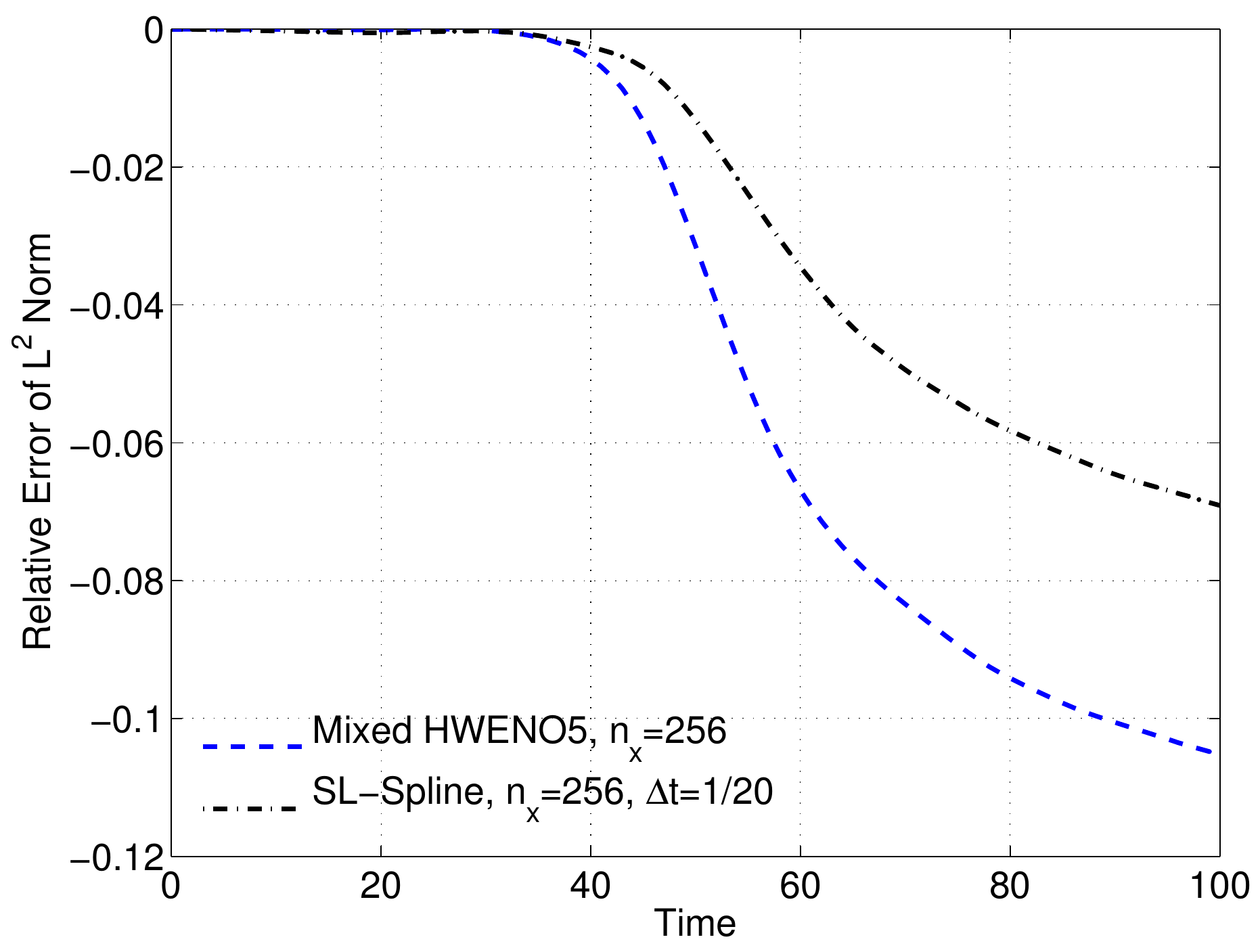}     & 
 \includegraphics[width=6.5cm]{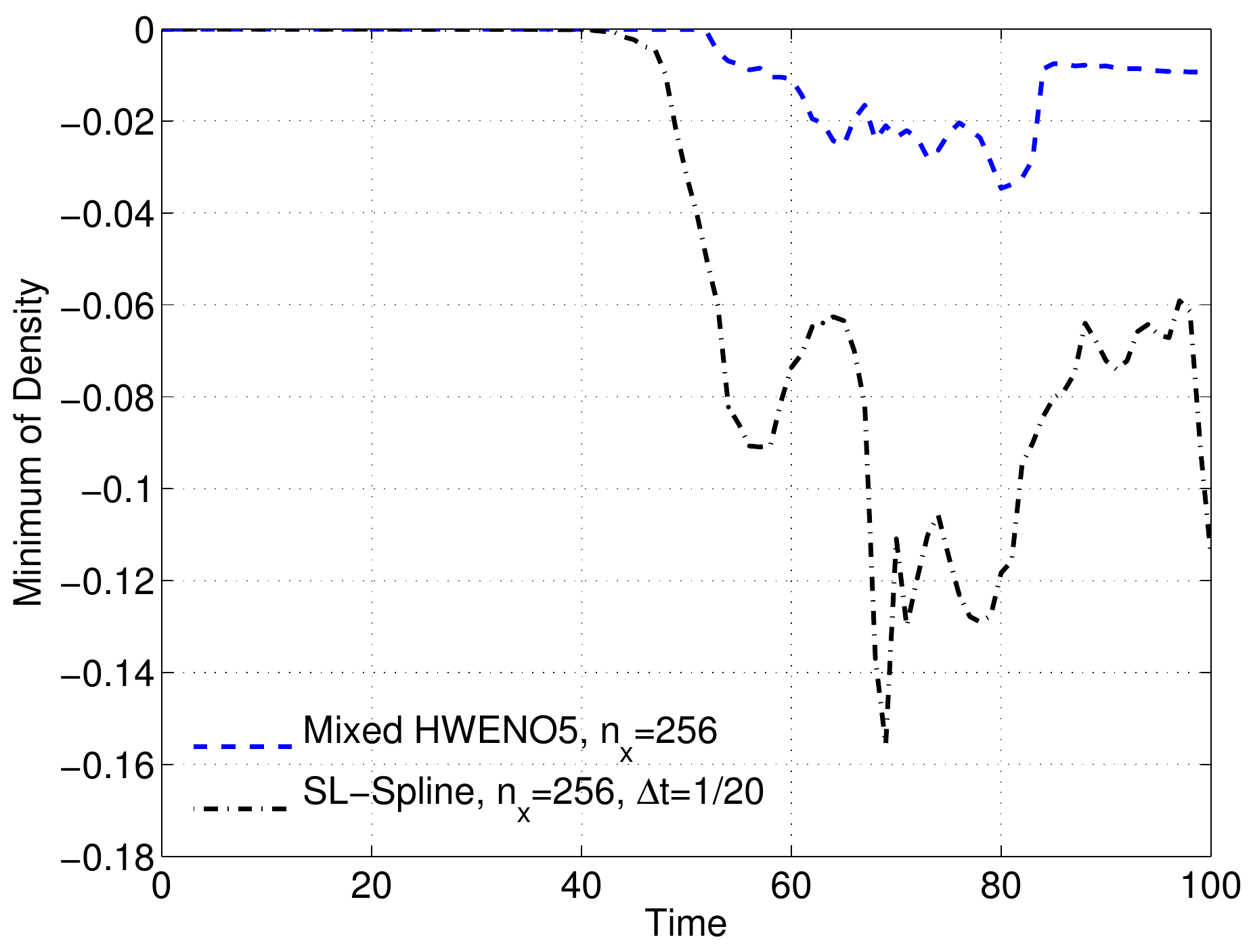}  
\\
  (c)       Relative Error of  $L^2$ Norm                       &   
 (d)  Minimum of Density   
 \end{tabular}
\caption{\label{fig:Diocotron_hybrid}Guiding center model: {\it Comparison between semi-Lagrangian with cubic spline method and mixed semi-Lagrangian/finite difference Hermite WENO5  method}}
 \end{center}
\end{figure}

\begin{figure}
\begin{center}
 \begin{tabular}{cc}
  \includegraphics[width=6.5cm]{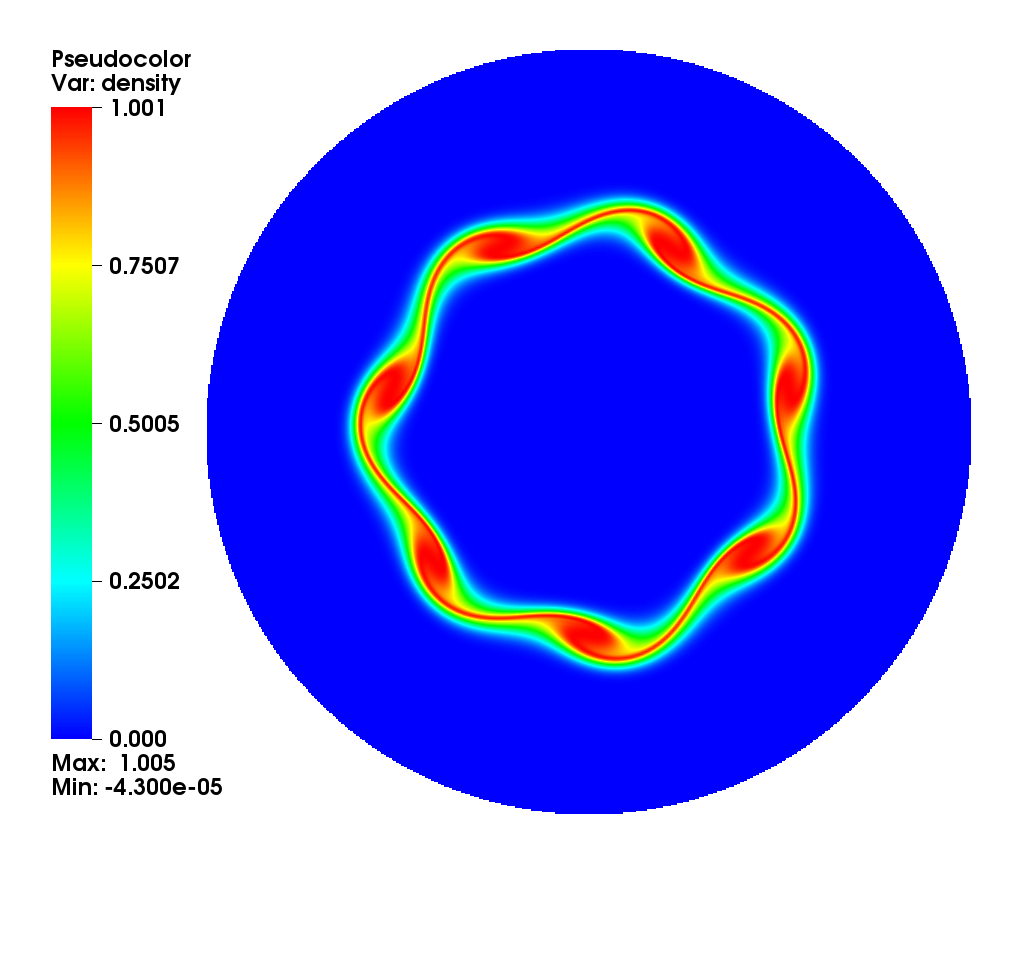} &   
 \includegraphics[width=6.5cm]{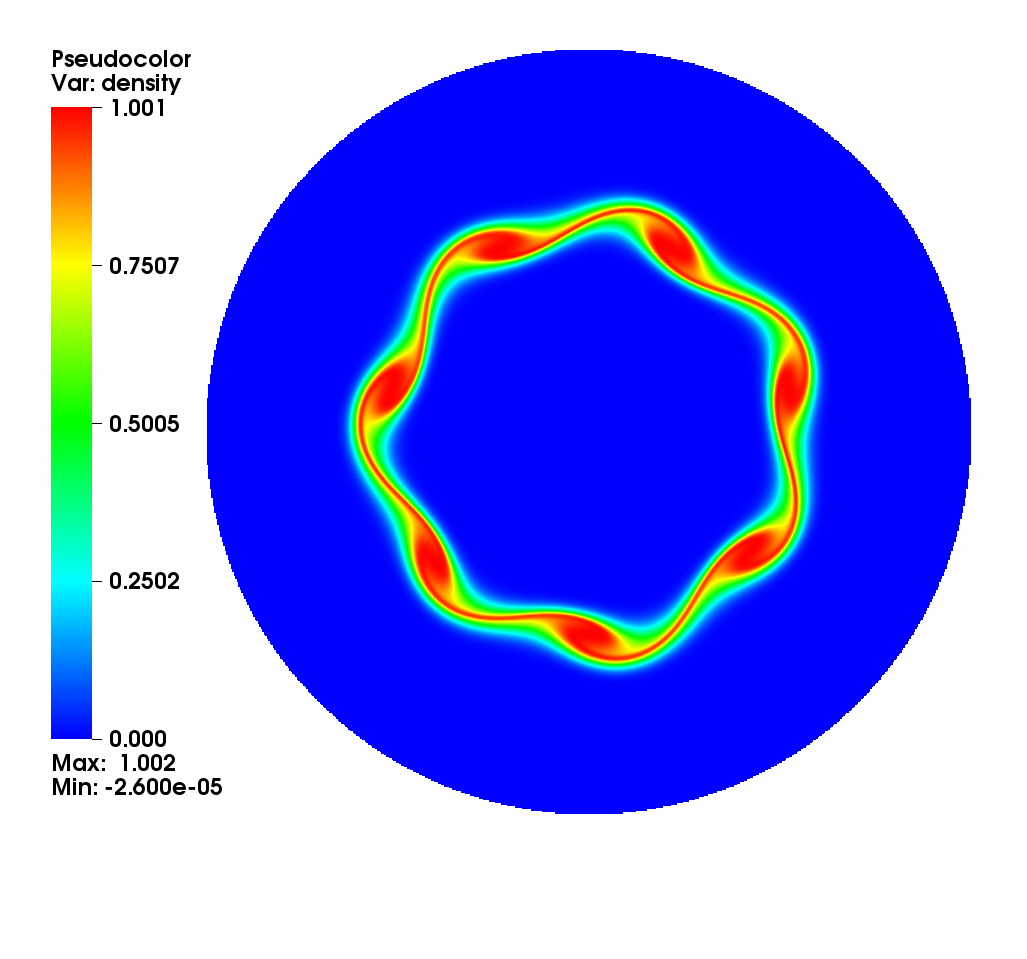} 
\\
%  (a) $t=40$                                             &    (b) $t=40$   \\
  \includegraphics[width=6.5cm]{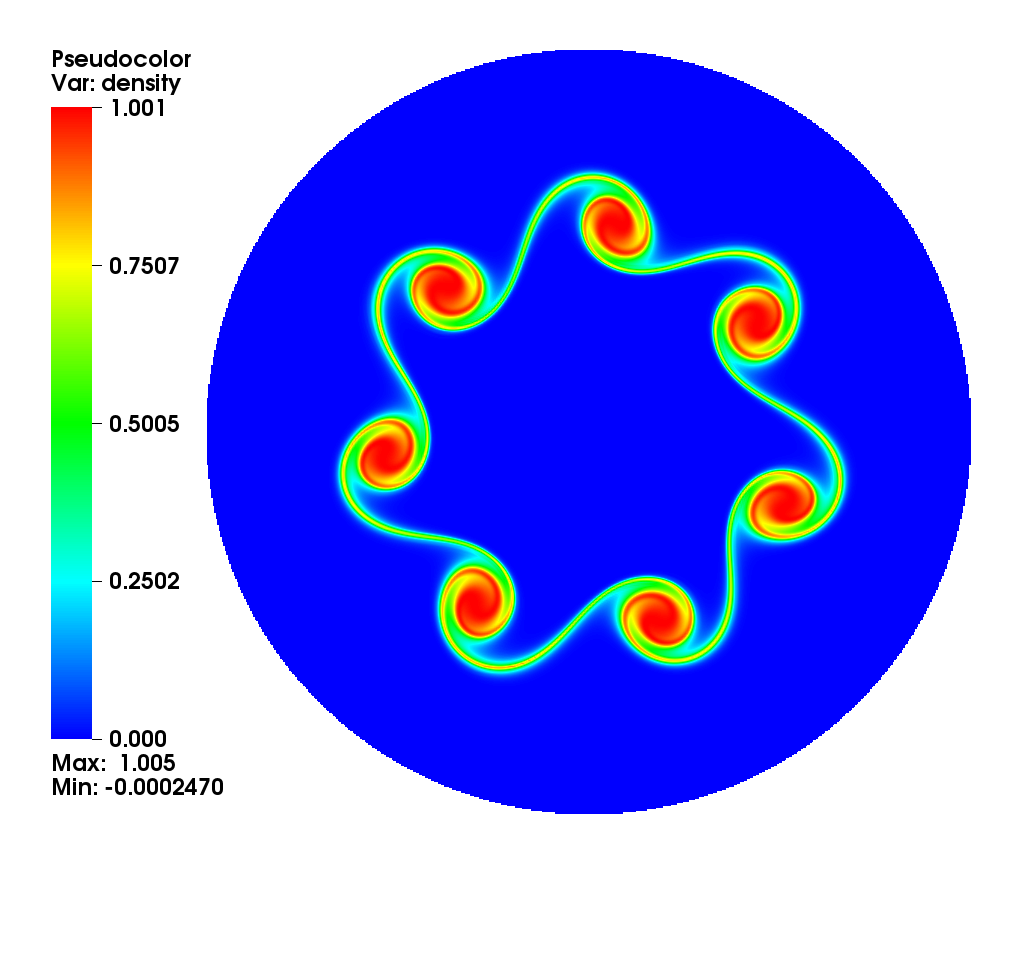} &    
\includegraphics[width=6.5cm]{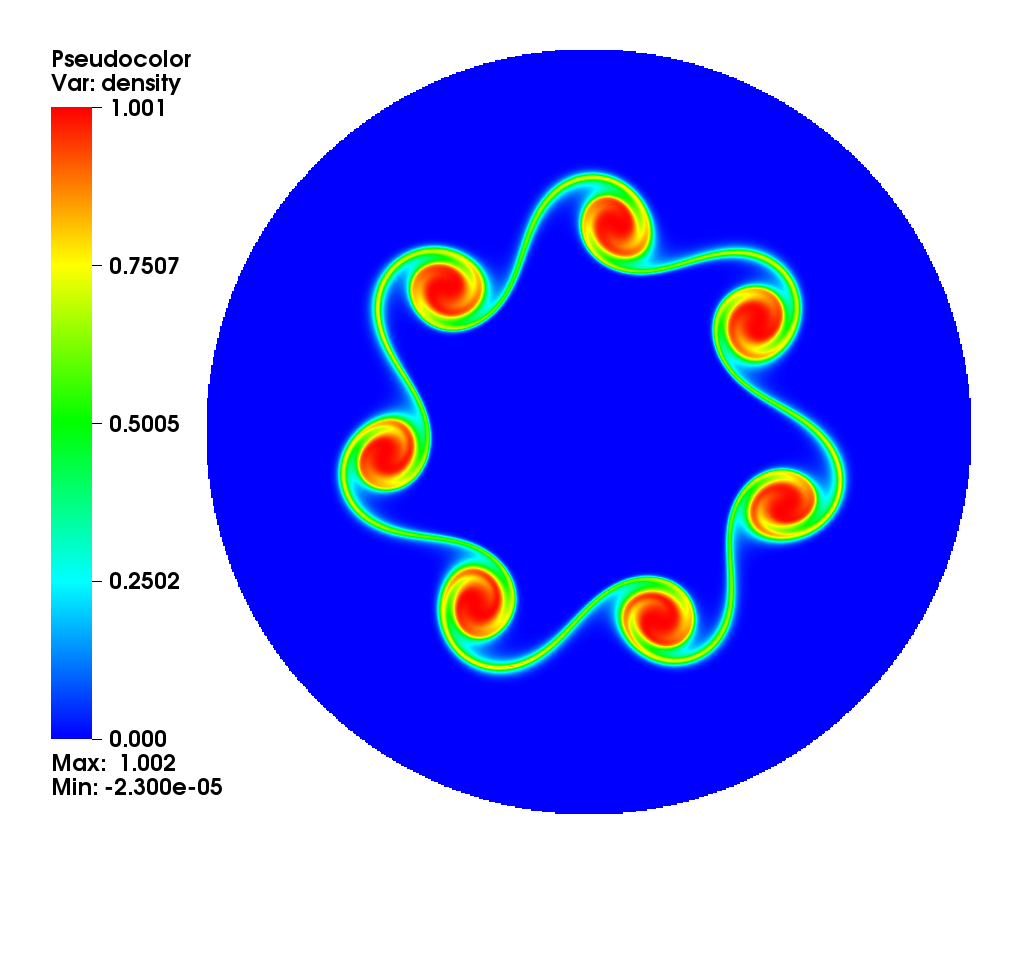} \\
 % (c)  $t=50$                                            &        (d) $t=50$   \\
  \includegraphics[width=6.5cm]{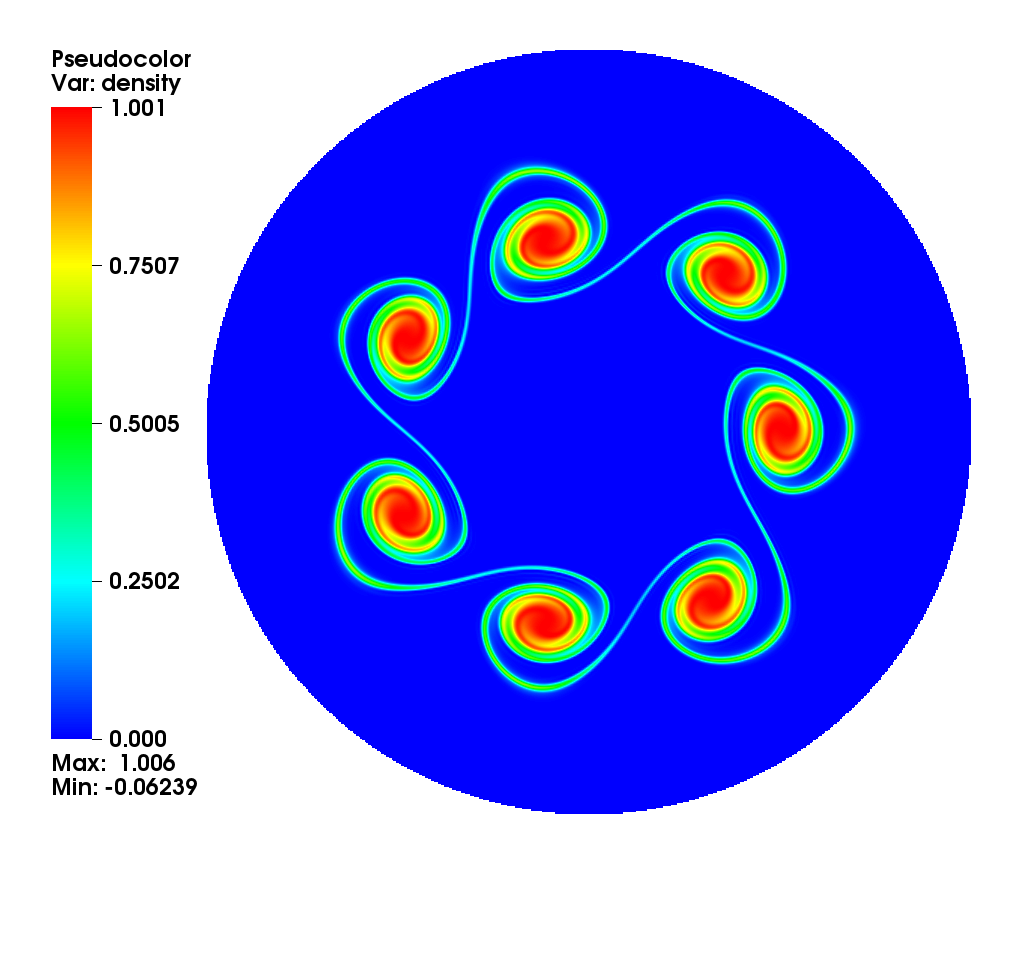} &    
\includegraphics[width=6.5cm]{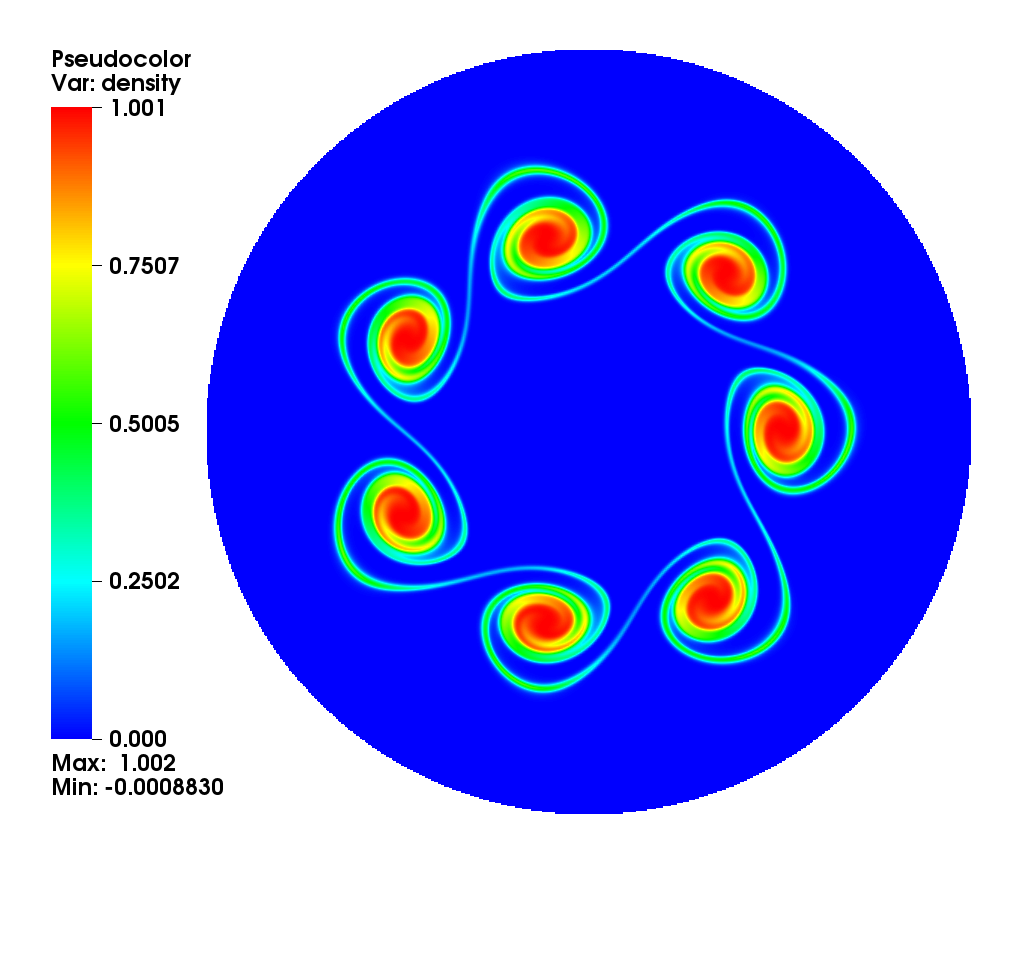} \\
%  (e)  $t=60$                                            &       (f)   $t=60$ 
(a) & (b) 
\end{tabular}
\caption{\label{fig:instability_GC}Guiding center model : {\it (a) semi-Lagrangian with cubic spline (b) mixed semi-Lagrangian/finite difference method with Hermite WENO5 at time $t= 40$, $50$ and $60$.}}
 \end{center}
\end{figure}

A comparison  between the semi-Lagrangian with cubic spline method and the mixed method for the diocotron instability is presented in Figures~\ref{fig:Diocotron_hybrid} and \ref{fig:instability_GC}. For a fair comparison an adaptive time step is also applied to the semi-Lagrangian scheme. We first choose a CFL number $\lambda\approx 2 $ during the linear phase, while  we take $\lambda\approx 0.5$  for the nonlinear phase and the number of points in space is $n_x=n_y=256$.

Using the semi-Lagrangian with cubic spline method, the relative error of mass oscillates a lot during the nonlinear phase, while relative error for the mixed method is more stable. We observe a very similar phenomenon for energy conservation. Unfortunately, relative error of $L^2$ norm for the mixed method is larger than that for the  semi-Lagrangian with cubic spline method, but the price to pay is to generate a non negligible negative values. Indeed, the WENO reconstruction  allow to control spurious oscillations generate from the discontinuous initial data and small structures.

Finally, the evolution of the density $\rho$ is presented in Figure~\ref{fig:instability_GC}. At first glance, we see the density of these two methods are very similar. At time $t=40$, small filaments appear and then  seven vortices are formed and move. 

Looking more carefully, we observe that  the numerical results obtained from the mixed method  is a little bit more dissipative than the ones obtained from the semi-Lagrangian methods with cubic spline interpolation {\it i.e.} small structures of density are more thin. However, the semi-Lagrangian method is  much more oscillatory than the mixed method, which can be observed from the minimum or maximum of density.

As a conclusion, although semi-Lagrangian method is less dissipative than  the mixed HWENO method, it involves too much numerical instabilities in nonlinear phase.
Therefore, the mixed method controlling spurious oscillations is more appropriate for long time simulation in plasma physics.

%%%%%%%%%%%%%%%%%%%%%%%%%%%%%%%%%%%%%%%%%%%%%%%%%%%%%%%%%%%%%%%%%%%%%%%
%%                                                                   %%
%%                                                                   %%
%%%%%%%%%%%%%%%%%%%%%%%%%%%%%%%%%%%%%%%%%%%%%%%%%%%%%%%%%%%%%%%%%%%%%%%
\section{Conclusion and perspective}
\label{sec:conc}
\setcounter{equation}{0}

In this paper, we have first developed a Hermite weighted essentially non-oscillatory reconstruction for semi-Lagrangian method and finite difference method respectively.

We illustrate that such a reconstruction is less dissipative than usual  weighted essentially non-oscillatory reconstruction.  Then we have compared our approach with the usual   semi-Lagrangian method with cubic spline and finite difference  WENO reconstruction.
The  semi-Lagrangian method is efficient and accurate for linear phase even with a large time step, however, it becomes less accurate for nonlinear phase and may lead to the wrong solution in some cases, for instance, the Beam test~\cite{bibCLM}.

The finite difference method is stable under the classical CFL condition, but it is much more stable in nonlinear phase and it conserves mass. We thus apply a mixed method using the  semi-Lagrangian method in linear phase and finite difference method during  the nonlinear phase, called mixed HWENO5 method.

We finally apply the mixed HWENO5 method to the simulation of the diocotron instability and observe that although the mixed HWENO5 method is a little more dissipative than the semi-Lagrangian with cubic spline method, but it is much more stable during the nonlinear phase.

The next step is now to apply our mixed method to more realistic and high dimensional plasma turbulence simulations, for instance, 4D Drift-Kinetic simulation~\cite{bibGV} or 5D Gyrokinetic simulation\cite{ref2}.

%%%%%%%%%%%%%%%%%%%%%%%%%%%%%%%%%%%%%%%%%%%%%%%%%%%%%%%%%%%%%%%%%%%%%%%
%%                                                                   %%
%%                                                                   %%
%%%%%%%%%%%%%%%%%%%%%%%%%%%%%%%%%%%%%%%%%%%%%%%%%%%%%%%%%%%%%%%%%%%%%%%
\section*{Acknowledgment} 
The authors are  partially supported by the European Research Council ERC Starting Grant 2009,  project 239983-\textit{NuSiKiMo}.
  
\bibliographystyle{plain}

\begin{flushleft} 
\signff 
\end{flushleft}
\vspace{-4.25cm}
\begin{flushright} 
\signcy 
\end{flushright}

\end{document}